\documentclass[11pt]{article}\textwidth 160mm\textheight 235mm
\usepackage{amsfonts}
\usepackage{amssymb}
\usepackage{graphicx}
\usepackage{amsmath}
\usepackage{amsthm}
\usepackage{enumitem}
\usepackage{tikz}
\usepackage{cancel}
\usepackage{todonotes}
\usepackage[normalem]{ulem}
\usetikzlibrary{matrix}
\usetikzlibrary{arrows,shapes}
\usepackage{cite}
\usepackage{hyperref}
\hypersetup{colorlinks=true, urlcolor=blue}
\expandafter\let\expandafter\oldproof\csname\string\proof\endcsname
\let\oldendproof\endproof
\renewenvironment{proof}[1][\proofname]{\oldproof[\ttfamily\scshape\bf #1.]
}{\oldendproof}

\def\O{{\Omega}}

\def\B{\mathbb{B}}
\def\R{{\rm I\!R}}
\def\N{{\rm I\!N}}
\def\ox{\bar{x}}
\def\oy{\bar{y}}
\def\oz{\bar{z}}
\def\ov{\bar{v}}
\def\ow{\bar{w}}
\def\ou{\bar{u}}

\def\ss{\scriptsize }
\def\ve{\varepsilon}
\def\epsilon{\varepsilon}

\def\hat{\widehat}

\def\emp{\emptyset}

\def\d{{\rm d}}
\def\sub{\partial}

\def\disp{\displaystyle}

\def\tto{\rightrightarrows}
\def\Hat{\widehat}

\def\Bar{\overline}
\def\ra{\rangle}
\def\la{\langle}
\def\ve{\varepsilon}

\def\co{\mbox{\rm co}\,}

\def\int{\mbox{\rm int}\,}
\def\gph{\mbox{\rm gph}\,}
\def\epi{\mbox{\rm epi}\,}

\def\dom{\mbox{\rm dom}\,}

\def\ker{\mbox{\rm ker}\,}

\def\dn{\downarrow}

\def\ph{\varphi}
\def\emp{\emptyset}
\def\st{\stackrel}
\def\oR{\Bar{\R}}
\def\lm{\lambda}
\def\olm{\bar\lambda}

\def\gg{\gamma}
\def\dd{\delta}
\def\al{\alpha}
\def\Th{\Theta}

\def\vt{\vartheta}
\def\sce{\setcounter{equation}{0}}
\def\ss{\scriptsize}
\def\vep{\varepsilon}

\def\fsub{\hat{\partial}}
\def\olambda{\bar{\lambda}}
\def\sm{\hbox{${1\over 2}$}}
\begin{document}
\newtheorem{Theorem}{Theorem}[section]
\newtheorem{Proposition}[Theorem]{Proposition}
\newtheorem{Remark}[Theorem]{Remark}
\newtheorem{Lemma}[Theorem]{Lemma}
\newtheorem{Corollary}[Theorem]{Corollary}
\newtheorem{Definition}[Theorem]{Definition}
\newtheorem{Example}[Theorem]{Example}
\renewcommand{\theequation}{{\thesection}.\arabic{equation}}
\renewcommand{\thefootnote}{\fnsymbol{footnote}}
\begin{center}
{\bf\Large Variational Analysis of Composite Models\\with Applications to Continuous Optimization}\\[2ex]
ASHKAN MOHAMMADI\footnote{Department of Mathematics, Wayne State University, Detroit, MI 48202, USA (ashkan.mohammadi@wayne.edu).}
BORIS S. MORDUKHOVICH\footnote{Department of Mathematics, Wayne State University, Detroit, MI 48202, USA (boris@math.wayne.edu).} and M. EBRAHIM SARABI\footnote{Department of Mathematics, Miami University, Oxford, OH 45065, USA (sarabim@miamioh.edu).}
\end{center}
\vspace*{-0.1in}
\small{\bf Abstract.} The paper is devoted to a comprehensive study of composite models in variational analysis and optimization the importance of which for numerous theoretical, algorithmic, and applied issues of operations research is difficult to overstate. The underlying theme of our study is a systematical replacement of conventional metric regularity and related requirements by much weaker metric subregulatity ones that lead us to significantly stronger and completely new results of first-order and second-order variational analysis and optimization. In this way we develop extended calculus rules for first-order and second-order generalized differential constructions with paying the main attention in second-order variational theory to the new and rather large class of fully subamenable compositions. Applications to optimization include deriving enhanced no-gap second-order optimality conditions in constrained composite models, complete characterizations of the uniqueness of Lagrange multipliers and strong metric subregularity of KKT systems in parametric optimization, etc.\vspace*{-0.15in}
\normalsize

\section{Introduction and Overview}\label{intro}
\sce\vspace*{-0.05in}

It has been well recognized in variational analysis, continuous optimization, and their various applications that composite models involving extended-real-valued functions of the type
\begin{equation}\label{CS}
(\vt\circ f)(x):=\vt\big(f(x)\big),\quad x\in\R^n,
\end{equation}
constitute a very convenient framework for developing both theoretical and algorithmic issues of constrained optimization with applications to practical modeling in operations research. Standard assumptions under which composite functions of type \eqref{CS} are investigated and applied in constrained optimization require that the mapping $f\colon\R^n\to\R^m$ be twice continuously differentiable (${\cal C}^2$-smooth), that the extended-real-valued function $\vt\colon\R^m\to\oR:=(-\infty,\infty]$ be lower semicontinuous (l.s.c.) and convex, and that the epigraphical set-valued mapping $H\colon\R^n\times\R\tto\R^m\times\R$ with $H(x,\al):=(f(x),\al)-\epi\vt$ be {\em metrically regular} around the point in question. We refer the reader to the book by Rockafellar and Wets \cite{rw} and the bibliographies therein for major facts on the theory and applications of such compositions known as (strongly) {\em amenable functions}, where the most perfect and complete results are obtained in the case of {\em fully amenable} compositions dealing with piecewise linear-quadratic outer functions $\vt$ in \eqref{CS}. The crucial metric regularity assumption mentioned above can be equivalently described as the {\em basic qualification condition} (or constraint qualification) expressed precisely at the reference point in question; see Section~\ref{1order} for more details and discussions.

It is important to emphasize that the possibility for $\vt$ to take the infinity value $\infty$ allows us to incorporate constraints in the composite unconstrained framework of \eqref{CS}. Indeed, while minimizing $\vt\circ f$ therein we automatically have the constraint $x\in\O:=\{x\in\R^n|\;f(x)\in\dom\vt\}$, where $\dom\vt:=\big\{y\in\R^m|\;\vt(y)<\infty\}$. On the other hand, the constrained optimization problem
\begin{equation*}
\mbox{minimize }\;\ph_{0}(x)\;\mbox{ subject to }\;g(x)\in\Th\subset\R^m
\end{equation*}
can be written in the unconstrained form \eqref{CS} via the functions $\vt(\alpha,y):=\alpha+\dd_\Th(y)$ and $f(x):=(\ph_{0}(x),g(x))$, where $\vt:\R^{m+1}\to\oR$ is l.s.c.\ and  $\dd_\Th$ is the indicator function of the set $\Th$ that equals to $0$ on $\Th$ and $\infty$ otherwise. This framework encompasses, in particular, problems of nonlinear programming (NLPs) where $\Th$ is a polyhedron and problems of conic programming where $\Th$ is a closed convex cone.

The aforementioned metric regularity and the equivalent notions of linear openness/covering and robust Lipschitzian behavior of set-valued mappings/multifunctions are largely investigated in variational analysis and are broadly applied to numerous topics in theoretical and computational optimization, equilibria, sensitivity analysis, optimal control. The reader may consult with the books by Borwein and Zhu \cite{bz}, Dontchev and Rockafellar \cite{dr}, Ioffe \cite{i}, Klatte and Kummer \cite{kk}, Mordukhovich \cite{m06,m18}, Penot \cite{pen}, and Rockafellar and Wets \cite{rw} together with the references and commentaries therein along with the enormous amount of other publications.\vspace*{0.05in}

Among the main intentions of this paper is to develop a new variational technique, which allows us to systematically replace metric regularity qualification conditions by much more subtle {\em metric subregularity} ones. The property of metric subregularity and its {\em calmness} equivalent for inverse mappings are less developed and applied in comparison with their robust metric regularity/Lipschitzian counterparts. While the latter properties admit complete characterizations, which open the door for developing comprehensive generalized differential calculus and various applications to optimization, stability, and other areas of nonlinear analysis and operations research, the study of metric subregularity and calmness is essentially more involved and the obtained results are by far more limited. On the other hand, such properties hold in many important settings where their robust counterparts are out of reach. Some results and discussions related to these topics can be found in the aforementioned monographs. We also refer the reader to the papers by Arag\'on Artacho and Geoffroy \cite{ag}, Chieu et al. \cite{chnt}, Drusvyatskiy et al. \cite{dmn}, Gfrerer \cite{g}, Gfrerer and Outrata \cite{go16}, Henrion and Outrata \cite{ho}, Ioffe and Outrata \cite{io}, Kruger \cite{k}, Li and Mordukhovich \cite{lm}, Ye and Ye \cite{yy}, and Zheng and Ng \cite{zn} with the additional bibliographies therein for various developments on metric subregularity and calmness properties and their applications to optimality and stability conditions, error bounds, and convergence of numerical algorithms.

The major goals of this paper are largely different from those considered in the literature on metric regularity, calmness, and their applications. Along with deriving refined calculus rules of first-order variational analysis under enhanced metric subregularity qualification conditions, we mainly concentrate on {\em second-order} variational analysis of compositions \eqref{CS} for a novel class of {\em fully subamenable} functions, where $\vt\colon\R^m\to\oR$ is l.s.c., convex, and piecewise linear-quadratic, where $f\colon\R^m\to\R^m$ is twice differentiable at the reference point, and where certain metric subregularity qualification conditions are satisfied. This class of functions is a direct extension of the fully amenable one, introduced by Rockafellar \cite{r88}, who imposed more restrictive metric regularity constraint qualifications, which allowed him to employ robust machinery of generalized differentiation; see the book \cite{rw} by Rockafellar and Wets for further developments and applications. The subregularity framework adopted here covers an essentially larger territory (including important settings where our subregularity conditions hold automatically), while it requires to develop fairly different techniques. The suggested variational approach allows us not only to significantly extend the known second-order calculus rules, but also to simplify their proofs and to obtain important results, which are new even in the case of metric regularity constraint qualifications.

Arguing in this way with a systematic usage of metric subregularity and variational ideas, we establish enhanced first-order chain rules for subderivatives and subdifferential mappings in rather general settings and then proceed with developing advanced second-order calculus for fully subamenable compositions. Our analysis in this paper mainly addresses second-order generalized differential constructions of the {\em primal} type (second and parabolic subderivatives) and {\em primal-dual} type (subgradient graphical derivatives) along with the associated second-order geometric objects. For all of them we develop extended second-order calculus rules under metric subregularity qualification conditions  and also efficiently compute these constructions for fully
subamenable compositions.

Besides employing variational techniques to derive the major second-order calculus and computational results, this paper develops a variety of applications to composite optimization and variational stability under the metric subregularity constraint qualification. Among such applications we mention no-gap second-order optimality conditions for composite optimization problems, characterizations of the uniqueness of Lagrange multipliers in the corresponding KKT (Karush-Kuhn-Tucker) systems, and characterizations of robust isolated calmness for solution maps to parameterized generalizes equations.\vspace*{0.05in}

The rest of the paper is organized as follows. In Section~\ref{prelim} we recall the basic notions of variational analysis that are broadly used below and present some preliminary results. Section~\ref{1order} is devoted to {\em first-order variational analysis} of compositions \eqref{CS}. We introduce here a new {\em metric subregularity qualification condition} for \eqref{CS}, discuss its relationships with the known ones in this setting, and employ it to derive enhanced first-order chain rules of the {\em equality} type for subderivatives and subgradient mappings.

In Section~\ref{subamen} we start developing {\em second-order variational analysis} of composite functions with the main emphasis on a novel class of {\em fully subamenable compositions}. First we define the critical cone of an l.s.c.\ function and prove that the critical cone agrees with the domain of the second subderivative under a certain second-order sufficient condition. Then we show that this second-order condition 
automatically holds if the function in question is fully subamenable.

In the Section~\ref{2subamen}, we prove that any fully subamenable function enjoys the powerful property of {\em twice epi-differentiability} and derive precise formulas for calculating its second subderivative. 
Section~\ref{2opt} develops applications of the obtained second-order calculus to deriving {\em no-gap second-order necessary and sufficient optimality conditions} for optimization problems with constraints given by fully subamenable compositions. The underlying metric subregularity qualification condition serves now as a refined constraint qualification ensuring in generality the validity of first-order optimality/stationarity conditions in the normal/KKT form and then the fulfillment of the aforementioned no-gap second-order optimality conditions in the case of fully subamenable constraints.

Section~\ref{gra-der} addresses yet another second-order generalized derivative for l.s.c.\ functions known as the {\em proto-derivative}. { For  fully subamenable compositions, we derive a precise calculation formula representing the proto-derivative of their subgradient mappings} entirely via the given data. 

Section~\ref{param} is devoted to applications of the above computations of the proto-derivative to the important and somewhat interconnected issues of {\em parametric composite optimization} concerning the {\em uniqueness of Lagrange multiplies} and {\em strong metric subregularity} of solution maps to KKT systems. In this way we derive, in particular, constructive second-order characterizations of strong metric subregularity of solution maps to the KKT systems associated with \eqref{CS}, where the outer function $\vt$ is piecewise linear-quadratic.  One of the these characterizations provides a new and much simpler proof of the recent result due to Burke and Engle \cite{be}. Concluding remarks in Section~\ref{conc} briefly summarize major results of the paper and discuss some directions of the future research.\vspace*{0.05in}

\color{black}Notation and terminology of this paper are standard in variational analysis and optimization. They are mainly taken, together with the preliminaries in Section~\ref{prelim}, from the books by Rockafellar and Wets \cite{rw} and by Mordukhovich \cite{m18}. For the reader's convenience and notational unification we usually use small Greek letters to denote scalar and extended-real-valued functions, small Latin letters for vectors and single-valued mappings/vector functions, and capital letters for sets, set-valued mappings, and matrices. Given a nonempty set $\O$ and vector $x$ in the Euclidean space $\R^n$, denote by ${\rm dist}(x;\O)$ the distance between $x$ and $\O$. The notation $\co\O$ stands for the convex hull of $\O$, while the symbol $x\st{\O}{\to}\ox$ indicates that $x\to\ox$ with $x\in\O$. By $\B$ we denote the closed unit ball in the space in question and by $\B_r(x):=x+r\B$ the closed ball centered at $x$ with radius $r>0$. The symbol $\{v\}^\bot$ with $v\in\R^n$ stands for the orthogonal subspace $\{w\in\R^n|\; \la w,v\ra=0\}$. The space of $n\times n$ matrices is denoted by $\R^{n\times n}$. As always, the vector quantity $x=o(t)$ with $t>0$ means that $\|x\|/t\to 0$ as $t\dn 0$. Recall also that $\R_+$ and $\R_-$ signify, respectively, the collection of nonnegative and nonpositive real numbers, and that $\N:=\{1,2,\ldots\}$. Given a scalar function $\ph\colon\R^n\to\R$, denote by $\nabla\ph(\ox)$ and $\nabla^2\ph(\ox)$ the gradient and Hessian of $\ph$ at $\ox$, respectively. If $f=(f_1,\ldots,f_m)\colon\R^n\to\R^m$ is a vector function twice differentiable at $\ox\in\R^n$, its second derivative at this point, denoted by $\nabla^2f(\ox)$, is a bilinear mapping from $\R^n\times\R^n$ into $\R^m$. In what follows we use the notation $\nabla^2f(\ox)(w,v)$, which means that
\begin{equation*}
\nabla^2f(\ox)(w,v)=\big(\big\la\nabla^2 f_1(\ox)w,v\big\ra,\ldots,\big\la\nabla^2 f_m(\ox)w,v\big\ra\big)\;\mbox{ for all }\;v,w\in\R^n.
\end{equation*}\vspace*{-0.3in}

\section{Basic Definitions and Preliminaries}\label{prelim}
\sce\vspace*{-0.05in}

We begin with recalling the well-known notions of variational analysis and generalized differentiation that are largely utilized and studied throughout the entire paper. Given a nonempty set $\O\subset\R^n$ with $\ox\in\O$, the (Bouligand-Severi) {\em tangent/contingent cone} $T_ \O(\ox)$ to $\O$ at $\ox\in\O$ is
\begin{equation}\label{tan1}
T_\O(\ox):=\big\{w\in\R^n\big|\;\exists\,t_k{\downarrow}0,\;\exists\,w_k\to w\;\mbox{ as }\;k\to\infty\;\;\mbox{with}\;\;\ox+t_kw_k\in\O\big\}.
\end{equation}
We say that a tangent vector $w\in T_\O(\ox)$ is {\em derivable} if there exists $\xi\colon[0,\ve]\to\O$ with $\ve>0$, $\xi(0)=\ox$, and $\xi'_+(0)=w$, where the right derivative $\xi'_+$ of $\xi$ at $0$ is defined by
$$
\xi'_+(0):=\lim_{t\dn 0}\frac{\xi(t)-\xi(0)}{t}.
$$
The set $\O$ is {\em geometrically derivable} at $\ox$ if every tangent vector $w$ to $\O$ at $\ox$ is derivable. This class of sets is sufficiently broad and includes, in particular, prox-regular sets and subdifferential graphs of convex piecewise linear-quadratic functions that are widely used in what follows.

The (Fr\'echet) {\em regular normal cone} to $\O$ at $\ox\in\O$ is
\begin{equation}\label{rn}
\Hat N_\O(\ox):=\Big\{v\in\R^n\Big|\;\limsup_{x\st{\O}{\to}\ox}\frac{\la v,x-\ox\ra}{\|x-\ox\|}\le 0\Big\},
\end{equation}
which can be equivalently described as the polar of the contingent cone $\Hat N_\O(\ox)=T_\O(\ox)^*$. The (Mordukhovich) {\em basic/limiting normal cone} to $\O$ at $\ox$ is defined by 

\begin{eqnarray}\label{ln}
N_\O(\ox)=\big\{v\in\R^n\;\big|\;\exists\,x\st{\O}{\to}\ox,\;v_k\to v\;\;\mbox{with}\;\;v_k\in\Hat N_\O(x_k)\;\mbox{ for all }\;k\in\N\big\}.
\end{eqnarray}
If the set $\O$ is convex, then both constructions \eqref{rn} and \eqref{ln} reduce to the classical normal cone of convex analysis. We say that $\O$ {\em normally regular} at $\ox\in\O$ if $\Hat N_\O(\ox)=N_\O(\ox)$.

Given a function $\ph\colon\R^n\to\oR$ with $\ox\in\dom\ph:=\{x\in\R^n|\;\ph(x)<\infty\}$, the {\em regular subdifferential} and the {\em limiting subdifferential} of $\ph$ at $\ox$ are defined via the regular \eqref{rn} and limiting \eqref{ln} normal cones to the epigraph $\epi\ph:=\{(x,\al)\in\R^n\times\R|\;\al\ge\ph(x)\}$ of $\ph$ by
\begin{eqnarray}\label{sub}
\begin{array}{ll}
&\Hat\partial\ph(\ox):=\big\{v\in\R^n\big|\;(v,-1)\in\Hat N_{{\scriptsize\epi\ph}}\big(\ox,\ph(\ox)\big)\big\},\\
&\partial\ph(\ox):=\big\{v\in\R^n\big|\;(v,-1)\in N_{{\scriptsize\epi\ph}}\big(\ox,\ph(\ox)\big)\big\},
\end{array}
\end{eqnarray}
respectively. We say that $\ph$ is {\em lower regular} at $\ox$ if $\partial\ph(\ox)=\Hat\partial\ph(\ox)$.  Note that this notion does not reduce to the subdifferential regularity for non-Lipschitzian functions $\ph\colon\R^n\to\oR$. On the other hand, the normal regularity of sets agrees for locally closed sets in finite-dimensional spaces (while not in infinite dimensions) with the (directional) Clarke regularity; see Corollary~6.29 in the book by Rockafellar and Wets \cite{rw}. In what follows we prefer to use the regularity terminology from the books by Mordukhovich \cite{m06,m18}.

For a set-valued mapping $F\colon\R^n\rightrightarrows\R^m$ with its domain and graph given by
$$
\dom F:=\big\{x\in\R^n\big|\;F(x)\ne\emp\big\}\;\mbox{ and }\;\gph F:=\big\{(x,y)\in\R^n\times\R^m\big|\;y\in F(x)\big\},
$$
respectively, its {\em graphical derivative} $DF(\ox,\oy)\colon\R^n\tto\R^m$ at $(\ox,\oy)\in\gph F$ is defined via the tangent cone \eqref{tan1} to its graph at $(\ox,\oy)$ by
\begin{equation}\label{gder}
DF(\ox,\oy)(u):=\big\{v\in\R^m\big|\;(u,v)\in T_{\scriptsize{\gph F}}(\ox,\oy)\big\}\;\mbox{ for all }\;u\in\R^n,
\end{equation}
while its {\em coderivative} $D^*F(\ox,\oy)\colon\R^n\tto\R^m$  at this point is given by
\begin{equation}\label{cod}
D^*F(\ox,\oy)(v):=\big\{u\in\R^n\big|\;(u,-v)\in N_{\scriptsize{\gph F}}(\ox,\oy)\big\},\quad v\in\R^m.
\end{equation}

One of the central well-posedness concepts in nonlinear analysis with great many applications is the {\em metric regularity} of $F\colon\R^n\tto\R^m$ around $(\ox,\oy)\in\gph F$ postulated as the existence of a constant $\kappa\in\R_+$ and neighborhoods $U$ of $\ox$ and $V$ of $\oy$ such that
\begin{equation}\label{metreg}
{\rm dist}\big(x;F^{-1}(y)\big)\le\kappa\,{\rm dist}\big(y;F(x)\big)\quad\mbox{for all}\quad(x,y)\in U\times V.
\end{equation}
If $y=\oy$ in \eqref{metreg}, the mapping $F$ is said to be {\em metrically subregular} at $(\ox,\oy)$.

As mentioned in Section~\ref{intro}, metric regularity and the equivalent covering/linear openness and Lipschitzian properties of multifunctions admit complete characterizations via generalized differentiation. In this paper we employ the following result known as the {\em Mordukhovich/coderivative criterion}, which amounts to saying that a closed-graph mapping $F\colon\R^n\tto\R^m$ is metrically regular around $(\ox,\oy)\in\gph F$ with some modulus $\kappa\in\R_+$ if and only if
\begin{equation}\label{cod-cr}
{\rm ker}\,D^*F(\ox,\oy)=\big\{0\big\},
\end{equation}
where ${\rm ker}\,G:=\{v\in\R^m|\;0\in G(v)\}$ is the kernel of a set-valued mapping $G\colon\R^m\tto\R^n$; see Mordukhovich \cite{m93,m06,m18} and the book by Rockafellar and Wets \cite{rw} with the references therein for different proofs and discussions. The broad applicability of criterion \eqref{cod-cr} is largely due to {\em full calculus} available for the coderivative \eqref{cod} that is based  on {\em variational} and {\em extremal principles} of variational analysis.

Results of such a type are not available for metric subregularity in general. This  makes the justification of this property  significantly more challenging for major classes of constrained and composite optimization problems. Yet the latter property is satisfied in many important situations where metric regularity fails. This can happen, in particular,  for problems with a {\em polyhedral} structure including the composite function \eqref{CS} with a convex piecewise linear-quadratic outer function $\vt$ and an affine inner function $f$.

We say that $\ph\colon\R^n\to\oR$ is {\em piecewise linear-quadratic} (PWLQ) if $\dom\ph=\cup_{i=1}^{s}\O_i$ with $\O_i$ being polyhedral convex sets for $i=1,\ldots,s$, and if $\ph$ has a representation of the form
\begin{equation}\label{PWLQ}
\ph(x)=\sm\langle A_ix,x\rangle+\langle a_i,x\rangle+\alpha_i\quad\mbox{for all}\quad x \in\O_i,
\end{equation}
where $A_i$ is an $n\times n$ symmetric matrix, $a_i\in\R^n$, and $\alpha_i\in\R$ for all $i=1,\ldots,s$. Recall also that $\ph\colon\R^n\to\oR$ is (locally) {\em Lipschitz continuous} around $\ox\in\dom\ph$ {\em relative} to some set $\O\subset\dom\ph$ if there exist a constant $\ell\in\R_+$ and a neighborhood $U$ of $\ox$  such that
\begin{equation*}
|\ph(x)-\ph(u)|\le\ell\|x-u\|\quad\mbox{for all }\quad x,u\in\O\cap U.
\end{equation*}
Piecewise linear-quadratic and indicator functions are simple albeit important examples of extended-real-valued functions that are Lipschitz continuous relative to their domains around every point therein.

A more delicate property (than the local Lipschitz continuity) of extended-real-valued functions $\ph\colon\R^n\to\oR$ {\em at} the point in question is the {\em calmness from below} of $\ph$ at $\ox\in\dom\ph$ meaning that there exist a constant $\ell\in\R_+$ and a neighborhood $U$ of $\ox$ such that
\begin{equation}\label{calmfrombelow}
\ph(x)\ge\ph(\ox)-\ell\|x-\ox\|\quad\mbox{for all}\quad x\in U.
\end{equation}

The following proposition describes useful consequences of this notion that is employed below.\vspace*{-0.05in}

\begin{Proposition}[\bf existence of subgradients]\label{nes} Let $\ph\colon\R^n\to\oR$ be l.s.c.\ around $\ox\in\dom\ph$ and calm from below at this point with constant $\ell\in\R_+$. Then we have $\partial\ph(\ox)\cap\ell\B\ne\emp$. This ensures that $\sub\ph(\ox)\ne\emp$ for any function $\ph$  that is piecewise linear-quadratic.
\end{Proposition}\vspace*{-0.15in}
\begin{proof} The calmness from below \eqref{calmfrombelow} clearly implies that the function
\begin{equation*}
\psi(x):=\ph(x)+\ell\|x-\ox\|,\quad x\in\R^n,
\end{equation*}
attains its local minimum at $\ox$. Then the subdifferential Fermat rule tells us that $0\in\partial\psi(\ox)$. Taking into account that $\ph$ is l.s.c.\ around $\ox$ while the function $x\mapsto\|x-\ox\|$ is obviously locally Lipschitzian, we use the semi-Lipschitzian sum rule from Theorem~2.33(c) in Mordukhovich \cite{m06} to conclude that $0\in\partial\psi(\ox)\subset\partial\ph(\ox)+\ell\B$, which verifies the first assertion of the proposition.

If $\ph$ is piecewise linear-quadratic, it follows from \eqref{PWLQ} that
\begin{equation*}
\ph(x)=\min_{1\le i\le s}\big\{\ph_{i}(x)+\dd_{\O_i}(x)\big\}\;\mbox{ with }\ph_{i}(x):=\sm\la A_ix,x\ra+\la a_i,x\ra+\al_i
\end{equation*}
and with $\O_{i}$ taken from \eqref{PWLQ}. This representation easily implies that $\ph$ is calm from below with some constant $\kappa\in\R_+$ at any point of its domain. This yields  $\sub\ph(\ox)\ne\emp$ for any  $\ox\in\dom\ph$ due to the  first assertion of this proposition.
\end{proof}\vspace*{-0.05in}

Note that the first part of Proposition~\ref{nes} can be deduced from Proposition~8.32 in Rockafellar and Wets \cite{rw} by using a different approach, while the second part is verified in Proposition~10.21 of that book when $\ph$ is assumed in addition to be convex.\vspace*{0.05in}

To finish with the first-order constructions, recall that the {\em subderivative} of $\ph\colon\R^n\to\oR$ at $\ox\in\dom\ph$ is a positively homogeneous function ${\rm d}\ph(\ox)\colon\R^n\to[-\infty,\infty]$ defined by
\begin{equation}\label{subder}
{\mathrm d}\ph(\ox)(\ow):=\disp\liminf_{\substack{t\dn 0\\w\to\ow}}{\frac{\ph(\ox+tw)-\ph(\ox)}{t}}\;\mbox{ for any }\;\ow\in\R^n.
\end{equation}
There is the well-known duality correspondence between the subderivative \eqref{subder} and the regular subdifferential of $\ph$ at $\ox$ taken from \eqref{sub}:
\begin{equation}\label{dual}
\Hat\partial\ph(\ox)=\big\{v\in\R^n\big|\;\la v,w\ra\le\d\ph(\ox)(w),\quad w\in\R^n\big\}.
\end{equation}
If $\ph$ is convex and piecewise linear-quadratic, then for any $\ox\in\dom\ph$ we get from Proposition~10.21 in Rockafellar and Wets \cite{rw} that
\begin{equation}\label{dfs}
\dom\d\ph(\ox)=T_{\ss\dom\ph}(\ox)=\bigcup_{i\in I(\ox)}T_{\O_i}(\ox)\quad\mbox{with}\quad I(\ox):=\big\{i\in\{1,\ldots,s\}\big|\;\ox\in\O_i\big\}.
\end{equation}
Furthermore, for any $w\in\dom\d\ph(\ox)$ there exists an index $i\in I(\ox)$ such that $w\in T_{\O_i}(\ox)$ and
\begin{equation}\label{dfs2}
\d\ph(\ox)(w)=\la A_i\ox+a_i,w\ra.
\end{equation}

Proceed now with the primal second-order constructions studied in this paper. Form the parametric family of second-order difference quotients of $\ph$ at $\ox\in\dom\ph$ for some $\ov\in\R^n$ by
\begin{equation*}
\Delta_t^2\ph(\ox,\ov)(w):=\dfrac{\ph(\ox+tw)-\ph(\ox)-t\langle\ov,w\rangle}{\sm t^2}\quad\mbox{with }\;w\in\R^{n},\;t>0.
\end{equation*}
The {\em second-order subderivative} $\d^2\ph(\ox,\ov)\colon\R^n\to[-\infty,\infty]$ of $\ph$ at $\ox$ for $\ov$ is defined by
\begin{equation}\label{ssd}
\d^2\ph(\bar x,\ov)(w):=\liminf_{\substack{t\dn 0\\u\to w}}\Delta_t^2 \ph(\ox,\ov)(u),\quad w\in\R^{n}.
\end{equation}
It is said that $\ph\colon\R^n\to\oR$ is {\em twice epi-differentiable} at $\bar x$ for $\ov$ if the second-order difference quotients $\Delta_t^2\ph(\ox,\ov)$ epi-converge to $\d^2\ph(\bar x,\ov)$ as $t\downarrow 0$; see Definition~7.1 in Rockafellar and Wets \cite{rw}. If in addition the second subderivative \eqref{ssd} is a proper function, then $\ph$ is called {\em properly twice epi-differentiable} at $\bar x$ for $\ov$. Recall that the above properness means that $\d^2\ph(\bar x,\ov)(w)>-\infty$ for all $w\in\R^n$ with $\dom\d^2\ph(\bar x,\ov)\ne\emp$. Recall also from Proposition~7.2 in Rockafellar and Wets \cite{rw} that the twice epi-differentiability of $\ph$ at $\bar x$ for $\ov$  can be equivalently described as follows: for every $w\in\R^n$ and every sequence $t_k\downarrow 0$ there exists a sequence $w_k\to w$ such that
\begin{equation}\label{df02}
\Delta_{t_k}^2\ph(\bar x,\ov)(w_k)\to\d^2\ph(\bar x,\ov)(w)\;\mbox{ as }\;k\to\infty.
\end{equation}

Turning to second-order variational geometry and given $\O\subset\R^n$ with $\ox\in\O$, define the {\em second-order tangent set} to $\O$ at $\ox$ for a tangent vector $w\in T_\O(\ox)$  by
\begin{equation}\label{2tan}
T_\O^2(\ox,w)=\big\{u\in\R^n\big|\;\exists\,t_k{\downarrow}0,\;\exists\,u_k\to u\;\mbox{ as }\;k\to\infty\;\;\mbox{with}\;\;\ox+t_kw+\sm t_k^2u_k\in\O\big\}.
\end{equation}
A set $\O$ is said to be {\em parabolically derivable} at $\ox$ for $w$ if $T_\O^2(\ox,w)\ne\emp$ and for each $u\in T_\O^2(\ox,w)$ there exists $\xi\colon[0,\ve]\to\O$ with $\ve>0$, $\xi(0)=\ox$, and $\xi'_+(0)=w$ such that $\xi''_+(0)=u$, where
$$
\xi''_+(0):=\lim_{t\dn 0}\frac{\xi(t)-\xi(0)-t\xi'_+(0)}{\sm t^2}.
$$

We conclude this section with the following simple and useful fact about the second-order tangent set
to the domain of a convex piecewise linear-quadratic function.\vspace*{-0.05in}

\begin{Proposition}[\bf second-order tangent sets to domains of convex PWLQ functions]\label{sotsc}
Let $\ph\colon\R^n\to\oR$ be a convex piecewise linear-quadratic function with $\ox\in\dom\ph$, and let
$w\in\dom\d\ph(\ox)$. Then we have the representation
\begin{equation}\label{union}
T_{\ss\dom\ph}^2(\ox,w)=\bigcup_{i\in J(\ox,w)}T_{\O_i}^2(\ox,w),
\end{equation}
where $\O_i$ are taken from \eqref{PWLQ}, and where
\begin{equation}\label{union2}
J(\ox,w):=\big\{i\in I(\ox)\big|\;w\in T_{\O_i}(\ox)\big\}
\end{equation}
with the index set $I(\ox)$ defined in \eqref{dfs}.
\end{Proposition}\vspace*{-0.15in}
\begin{proof} The inclusion ``$\supset$" in \eqref{union} is an immediate consequence of the fact that $\O_i\subset\dom\ph$ for all $i=1,\ldots,s$.
The opposite inclusion therein follows from the representation of $\dom\ph$ as the finite union of the polyhedral convex sets $\O_i$.
\end{proof}\vspace*{-0.3in}

\section{First-Order Chain Rules under Metric Subregularity}\label{1order}\sce

This section is mainly devoted to first-order variational analysis of general composite functions of type \eqref{CS}. Among other issues, it contains the introduction and study of the new metric subregularity qualification condition, which plays a crucial role in both first-order and second-order variational analysis conducted in this paper as well as in their applications to optimization.

We begin with the following proposition, which formulates major qualification conditions (including the new one) for deriving chain rules for compositions \eqref{CS} and then establishes relationships between them. Note that in applications to optimization, where compositions \eqref{CS} are used for modeling constraints, such conditions play a role of constraint qualifications and are often labeled in this way. On the other hand, we consider them as a tool of analysis, which is not applied only to constrained optimization. Observe that the qualification conditions formulated below specifically address compositions of type \eqref{CS} and are not defined for arbitrary set-valued or single-valued mappings.\vspace*{-0.05in}

\begin{Proposition}[\bf relationships between qualification conditions]\label{equi} Let $f\colon\R^n\to\R^m$ be a single-valued mapping differentiable at some point $\ox\in R^n$, and let $\vt\colon\R^m\to\oR$ be a proper extended-real-valued function continuous relative to its domain. Form the composition \eqref{CS}, assume that $\vt(f(\ox))$ is finite, and consider the following qualification conditions:

{\bf(i)} The set-valued mapping $H\colon\R^n\times\R\tto\R^m\times\R$ defined by $H(x,\alpha):=(f(x),\alpha)-\epi\vt$ is metrically regular
around the point $\big((\ox,\vt(f(\ox))),(0,0)\big)$.

{\bf(ii)} The set-valued mapping $H$ defined in {\rm(i)} is metrically subregular at $\big((\ox,\vt(f(\ox))),(0,0)\big)$.

{\bf(iii)} The set-valued mapping $G\colon\R^n\tto\R$ defined by $G(x):=f(x)-\dom\vt$ is metrically subregular at the point $(\bar x,0)$.\\[0.5ex]
Then we always have the implications {\rm(i)}$\implies${\rm(ii)}$\implies${\rm(iii)}.
\end{Proposition}\vspace*{-0.15in}
\begin{proof} Implication (i)$\implies$(ii) is obvious. To verify the second implication, suppose that (ii) holds and thus find a constant $\kappa\ge 0$ as well as the neighborhoods $U$ of $\ox$ and $V$ of $\vt(f(\ox))$ such that
\begin{equation}\label{msqc}
{\rm dist}\big((x,\alpha);\epi\ph)\le\kappa\,{\rm dist}\big((f(x),\alpha);\epi\vt)\quad\mbox{for all}\quad(x,\alpha)\in U\times V.
\end{equation}
Picking $x\in U$ and $\ve>0$, we get a vector $y\in\dom\vt$ such that
$$
\|f(x)-y\|<{\rm dist}\big(f(x);\dom\vt\big)+\ve,
$$
which implies in turn the upper estimates
$$
\|y-f(\ox)\|\le\|f(x)-f(\ox)\|+\|y-f(x)\|\le 2\|f(x)-f(\ox)\|+\ve.
$$
Since $\vt$ is continuous at $f(\ox)$ relative to its domain, suppose by shrinking the neighborhood $U$ if necessary that
$\vt(y)\in V$. Using this together with \eqref{msqc} ensures the inequalities
\begin{eqnarray*}
{\rm dist}\big((x,\vt(y));\epi\ph\big)&\le&\kappa\,{\rm dist}\big((f(x);\vt(y));\epi\vt\big)\\
&\le&\kappa\|f(x)-y\|+|\vt(y)-\vt(y)|<\kappa\,{\rm dist}\big(f(x);\dom\vt\big)+\ve.
\end{eqnarray*}
On the other hand, for all $(x,\al)\in\R^n\times\R$ we always have
$$
{\rm dist}(x;\dom\ph)\le{\rm dist}\big((x,\alpha);\epi\ph\big).
$$
Combining the above inequalities with $G^{-1}(0)=\dom\ph$ and \eqref{msqc} brings us to the inequalities
$$
{\rm dist}\big(x;G^{-1}(0)\big)={\rm dist}(x;\dom\ph)\le\kappa\,{\rm dist}\big((x,\vt(y));\epi\ph\big)\le\kappa^2\,{\rm dist}\big(f(x);\dom\vt\big)+\kappa\ve,
$$
and therefore results in the distance estimate
$$
{\rm dist}\big(x;G^{-1}(0)\big)\le\kappa^2\,{\rm dist}\big(f(x);\dom\vt\big)\quad\mbox{for all}\quad x\in U.
$$
This clearly yields (iii) and hence completes the proof.
\end{proof}\vspace*{-0.05in}

It easily follows from the coderivative criterion for metric regularity \eqref{cod-cr} applied to the set-valued mapping $H$ in Proposition~\ref{equi} that the qualification condition in (i) of that proposition can be equivalently written, in the case where $f$ is smooth and $\vt$ is l.s.c., as
\begin{equation}\label{bqc}
\partial^\infty\vt\big(f(\ox)\big)\cap\ker\nabla f(\ox)^*=\{0\},
\end{equation}
where $\partial^\infty\ph(\ox)$ stands for the {\em singular subdifferential} of $\ph\colon\R^n\to\oR$ at $\ox\in\dom\ph$ defined by
\begin{equation}\label{ss}
\partial^\infty\ph(\ox):=\big\{v\in\R^n\big|\;(v,0)\in N_{{\scriptsize\epi\ph}}\big(\ox,\ph(\ox)\big)\big\},
\end{equation}
and where the symbol $^*$ stands for the matrix transposition/adjoint operator. { The metric regularity qualification condition \eqref{ss} serves as the most advanced constraint qualification in first-order variational analysis and its various applications to constrained and composite optimization problems; see, e.g., the monographs by Mordukhovich \cite{m06,m18} and by Rockafellar and Wets \cite{rw} with the references therein. For particular classes of optimization problems, condition \eqref{bqc} reduces to the classical Mangasarian-Fromovitz constraint qualification in nonlinear programming, and to the Robinson constraint qualification in conic programming, where $\vt$ is the indicator function of a closed convex cone. Note that when $\vt$ is convex, the metric regularity qualification condition \eqref{bqc} can be equivalently represented in the form
\begin{equation}\label{gf07}
N_{\ss\dom\vt}\big(f(\ox)\big)\cap\ker\nabla f(\ox)^*=\{0\}.
\end{equation}}
It has been well recognized in variational analysis and documented in the aforementioned monographs that the metric regularity qualification condition and its implementations in \eqref{bqc} and \eqref{gf07} ensure the first-order subdifferential chain rule
\begin{equation*}
\sub(\vt\circ f)(\ox)\subset\nabla f(\ox)^*\sub\vt\big(f(\ox)\big),
\end{equation*}
which holds as equality provided that $\vt$ is lower regular at $\oy:=f(\ox)$. Ioffe and Outrata \cite{io} significantly improved this chain rule by replacing in its assumptions the metric regularity qualification condition \eqref{bqc} (labeled in \cite{io} as the ``standard Mordukhovich-Rockafellar subdifferential qualification condition") with the metric subregularity of the {\em epigraphical} mapping $H$ in Proposition~\ref{equi}(ii). However, this result does not cover, e.g., a particular setting of our interest in this paper for the composite function \eqref{CS}, where $\vt$ is a convex piecewise linear-quadratic function and  $f$ is an affine mapping. The chain rule of the equality type for the latter case, important in first-order and second-order variational analysis, holds automatically if the metric subregularity of the epigraphical mapping $H$ in Proposition~\ref{equi}(ii) is replaced by a weaker (and simpler) metric subregularity of the {\em domain} mapping $G(x)=f(x)-\dom \vt$ from condition (iii) of that proposition. Indeed, the metric subregularity of this domain mapping  results from the classical Hoffman lemma \cite{hof} for the aforementioned case of \eqref{CS}; see Corollary~\ref{chpwlq} below. \vspace*{0.05in}

The above discussion motivates the following definition that plays a central role in both first-order and second-order developments and their applications in this paper.\vspace*{-0.05in}

\begin{Definition}[\bf metric subregularity qualification condition]\label{MSCQ} Given $f\colon\R^n\to\R^m$ and $\vt\colon\R^m\to\oR$, we say that the composition $\vt\circ f$ satisfies the {\sc metric subregularity qualification condition} {\rm(}MSQC{\rm)} at $\bar x\in\dom f$ with constant $\kappa\in\R_+$ if the mapping $x\mapsto f(x)-\dom\vt$ is metrically subregular at $(\bar x,0)$ with this constant.
\end{Definition}\vspace*{-0.05in}

Taking into account the structure of the mapping $f-\dom\vt$ in Definition~\ref{MSCQ} and using \eqref{metreg} with the fixed vector $y=\oy=0$, observe that the introduced MSQC with a prescribed constant $\kappa\in\R_+$ can be equivalently described via the existence of a neighborhood $U$ of $\ox$ such that the distance estimate
\begin{equation}\label{mscq}
{\rm dist}(x;\O)\le\kappa\,{\rm dist}\big(f(x);\dom\vt\big)\quad\mbox{with}\quad\O:=\big\{x\in\R^n\big|\;f(x)\in\dom\vt\big\}
\end{equation}
is satisfied for all $x\in U$ with the same number $\kappa$ as in Definition~\ref{MSCQ}.

{As mentioned earlier, while metric regularity can be fully characterized via the coderivative, finding conditions under which metric subregularity holds is rather challenging. When, however,  the outer function $\vt$ is convex piecewise linear-quadratic that yields the polyhedrality of $\dom\vt$, Gfrerer \cite[Theorem~2.6]{g2} achieved a rather simple and verifiable second-order condition to ensure the validity of the MSCQ \eqref{mscq}. His second-order condition in the framework of Definition~\ref{MSCQ} demonstrates that if in addition $f$ is twice differentiable at $\ox$ and if for every
$w\in\R^n\setminus \{0\}$ with $\nabla f(\ox)w\in T_{\ss\dom \vt}(f(\ox))$ and every $\lm\in N_ {\ss\dom\vt}\big(f(\ox)\big)\cap\ker\nabla f(\ox)^*$ the implication
\begin{equation}\label{gfe}
\la \lm,\nabla^2 f(\ox)(w,w)\ra\ge 0\implies \lm=0
\end{equation}
holds, then the mapping $x\mapsto f(x)-\dom\vt$ is metrically subregular at $(\bar x,0)$.}\vspace*{0.05in}

{The following example provides two cases of the composite functions \eqref{CS} with $\vt$ being convex piecewise linear-quadratic such that the metric regularity qualification condition fails while the metric subregularity qualification condition is satisfied.\vspace*{-0.05in}

\begin{Example}[\bf failure of metric regularity for composite functions]\label{fmr}{\rm Consider the class of extended-real-valued functions $\vt\colon\R^m\to\oR$ defined by
\begin{equation}\label{enlp}
\vt(y):=\sup_{u\in Z}\big\{\la y,z\ra-\sm\la Bz,z\ra\big\},
\end{equation}
where $Z$ is a polyhedral convex set, and where $B$ is an $m\times m$ positive-semidefinite symmetric matrix. This type of penalty functions was introduced by Rockafellar \cite{r}, where he formulated an important class of composite optimization problems under the name of {\em extended nonlinear programming} (ENLP). It is not hard to see that $\vt$ is a convex piecewise linear-quadratic function.

{\bf(a)} Consider the function $\vt$ from \eqref{enlp} with $m=2$,
\begin{equation*}
Z:=\R^2,\;\mbox{ and }\;B:=\begin{pmatrix}
1&0\\
0&0
\end{pmatrix}
\end{equation*}
and define the constraint mapping $f\colon\R^2\to\R^2$ by $f(x_1,x_2):=(x_1-x_2,0)$. It is easy to check that $\dom\vt=\{(y_1,y_2)|\;y_2= 0\}$. For $\ox:=(0,0)\in\R^2$ we have
$$
N_{\ss\dom\vt}\big(f(\ox)\big)\cap\ker\nabla f(\ox)^*=\big\{(\lm_1,\lm_2)\in\R^2\big|\;\lm_1=0\big\},
$$
which shows that the metric regularity qualification condition \eqref{bqc} fails at $\ox$. On the other hand, MSQC \eqref{mscq} holds at $\ox$ since the mapping $x\mapsto f(\ox)-\dom \vt$ is metrically subregular at $(\ox,0)$. This follows from the aforementioned Hoffman lemma since $f$ is an affine mapping and $\dom\vt$ is a polyhedral convex set.

{\bf(b)} Consider the function $\vt$ from \eqref{enlp} with $m=3$, $Z:=\R_+^3$, and $B:=0\in\R^{3\times 3}$.
Define the constraint mapping $f:\R^3\to\R^3$ by $f=(f_1,f_2,f_3)$ with $x=(x_1,x_2,x_3)\in\R^3$ and
$$
f_1(x):=x_1-\sm x^2_2,\quad f_2(x):= x_1-\sm x^2_3,\quad f_3(x):=-x_1-\sm x^2_1-\sm x^2_2-\sm x^2_3.
$$
Letting $\ox:=(0,0,0)\in\R^3$, deduce from $\dom\vt=\R_-^3$ that
\begin{equation}\label{fmrlm}
N_ {\ss\dom\vt}\big(f(\ox)\big)\cap\ker\nabla f(\ox)^*=\big\{(\lm_1,\lm_2,\lm_3)\in\R_+^3\big|\;\lm_1+\lm_2-\lm_3=0\big\},
\end{equation}
which implies that the metric regularity qualification condition \eqref{bqc} fails at $\ox$. For any $\lm=(\lm_1,\lm_2,\lm_3)$ from \eqref{fmrlm} and any $w=(w_1,w_2,w_3)\in\R^3\setminus\{0\}$ with $\nabla f(\ox)w\in T_{\ss\dom \vt}(f(\ox))=\R^3_-$ we conclude from the conditions
$$
\la\lm,\nabla^2f(\ox)(w,w)\ra=-\lm_3w_1^2-(\lm_1+\lm_3)w_2^2-(\lm_2+\lm_3)w_3^2 \ge 0
$$
that $\lm=0$, which confirms that \eqref{gfe} is satisfied. This tells us that the constraint mapping $x\mapsto f(x)-\dom\vt $ is metrically subregular at $(\ox,0)$.}
\end{Example}}

We begin our first-order analysis with a new chain rule for the subderivative of the composite function \eqref{CS}. The following theorem is significantly different from the best results in this direction given in Theorem~10.6 of Rockafellar and Wets \cite{rw}. The main improvement is the replacement of the metric regularity qualification condition \eqref{bqc} therein with the much weaker MSQC \eqref{mscq}. Also, contrary to that book, we establish the subderivative chain rule as {\em equality} without any subdifferential regularity. Finally, the smoothness requirement on the inner mapping $f$ in \eqref{CS} is weakened to its merely (Fr\'echet) differentiability at the point in question. On the other hand, we impose the local Lipschitz continuity of the outer function $\vt$ {\em relative} to its domain, which is not assumed in Rockafellar and Wets \cite[Theorem~10.6]{rw}.\vspace*{-0.05in}

\begin{Theorem}[\bf subderivative chain rules as equalities under metric subregularity]\label{regchain}
Let $f\colon\R^n\to\R^m$ be differentiable at $\ox\in\R^n$, and let $\vt\colon\R^m\to\oR$ be Lipschitz continuous around $f(\ox)$ relative to its domain. If MSQC \eqref{mscq} is satisfied at $\ox$ with some constant $\kappa\in\R_+$, then we have
\begin{equation}\label{dchain}
\d(\vt\circ f)(\ox)(w)=\d\vt\big(f(\ox)\big)\big(\nabla f(\ox)w\big),\quad w\in\R^m.
\end{equation}
\end{Theorem}\vspace*{-0.15in}
\begin{proof} Pick any $\ow\in\R^n$ and deduce from the differentiability of $f$ at $\ox$ that $\nabla f(\ox)w+\frac{o(t\|w\|)}{t}\to\nabla f(\ox)\ow$ as $t\dn 0$ and $w\to\ow$. Based on this and \eqref{subder}, we get the relationships
\begin{eqnarray*}
\disp\d(\vt\circ f)(\ox)(\ow)&=&\liminf_{\substack{t\dn 0\\w\to\ow}}\frac{\vt\big(f(\ox+t w)\big)-\vt\big(f(\ox)\big)}{t}\nonumber\\
&=&\liminf_{\substack{t\dn 0\\w\to\ow}}\frac{\vt\big(f(\ox)+t\nabla f(\ox)w+o(t\|w\|)\big)-\vt\big(f(\ox)\big)}{t}\nonumber\\
&=&\liminf_{\substack{t\dn 0\\w\to\ow}}\frac{\vt\big(f(\ox)+t(\nabla f(\ox)w+\frac{o(t\|w\|)}{t})\big)-\vt\big(f(\ox)\big)}{t}\nonumber\\
&\ge&\d\vt\big(f(\ox)\big)\big(\nabla f(\ox)\ow\big)\;\mbox{ whenever }\;\ow\in\R^n,
\disp
\end{eqnarray*}
which verify the inequality ``$\ge$" in \eqref{dchain}. Proceeding next with the proof of the opposite inequality in \eqref{dchain}, take any $w\in\R^n$ and observe from the Lipschitz continuity of $\vt$ around $f(\ox)$ relative to its domain that $\d\vt(f(\ox))(\nabla f(\ox)w)>-\infty$. Since the latter inequality is obvious if $\d\vt(f(\ox))(\nabla f(\ox)w)=\infty$, we may assume that the value $\d\vt(f(\ox))(\nabla f(\ox)w)$ is finite. By definition,  there exist sequences $t_k\dn 0$ and $v_k\to\nabla f(\ox)w$ such that
\begin{equation}\label{subd0}
\d\vt\big(f(\ox)\big)\big(\nabla f(\ox)w\big)=\lim_{\substack{n\to\infty}}\frac{\vt\big(f(\ox)+t_k v_k\big)-\vt\big(f(\ox)\big)}{t_k}.
\end{equation}
Remembering that $\big|\d\vt(f(\ox))(\nabla f(\ox)w)\big|<\infty$, suppose without lost of generality that $f(\ox)+t_k v_k\in\dom\vt$ for all $k\in\N$.
Then the imposed MSQC \eqref{mscq} at $\ox$ yields
\begin{equation*}
{\rm dist}(\ox+t_kw;\O)\le\kappa\,{\rm dist}\big(f(\ox+t_k w);\dom\vt\big),\quad k\in\N,
\end{equation*}
which in turn brings us to the relationships
\begin{eqnarray*}
{\rm dist}\Big(w;\frac{\O-\ox}{t_k}\Big)&\le&\frac{\kappa}{t_k}\,{\rm dist}\big(f(\ox)+t_k\nabla f(\ox)w+o(t_k);\dom\vt\big)\nonumber \\
&\le&\frac{\kappa}{t_k}\,\big\|f(\ox)+t_k\nabla f(\ox)w+o(t_k)-f(\ox)-t_k v_k\big\|\nonumber\\
&=&\kappa\,\Big\|\nabla f(\ox)w-v_k+\frac{o(t_k)}{t_k}\Big\|\;\mbox{ for all }\;k\in\N.
\disp
\end{eqnarray*}
This allows us to find vectors $w_k\in\frac{\O-\ox}{t_k}$ satisfying
$$
\|w-w_k\|\le\kappa\,\Big\|\nabla f(\ox)w-v_k+\frac{o(t_k)}{t_k}\Big\|+\frac{1}{k}.
$$
Thus we have $\ox+t_k w_k\in\O$ for all $k$ and $w_k\to w$ as $k\to\infty$. Combining these with \eqref{subd0}, we arrive at
\begin{eqnarray*}
\d\vt\big(f(\ox)\big)\big(\nabla f(\ox)w\big)&=&\disp\lim_{k\to\infty}\Big[\frac{\vt\big(f(\ox+t_k w_k)\big)-\vt\big(f(\ox)\big)}{t_k}+\frac{\vt\big(f(\ox)+t_k v_k)\big)-\vt\big(f(\ox+t_kw_k)\big)}{t_k}\Big]\\
&\ge&\disp\liminf_{k\to\infty}\frac{\vt\big(f(\ox+t_kw_k)\big)-\vt\big(f(\ox)\big)}{t_k}-\ell\lim_{k\to\infty}\Big\|\frac{f(\ox+t_k w_k)-f(\ox)}{t_k}-v_k\Big\|\\
&\ge&\d(\vt\circ f)(\ox)(w)-\ell\lim_{k\to\infty}\Big\|\nabla f(\ox)w_k+\frac{o(t_k)}{t_k}-v_k\Big\|=\d(\vt\circ f)(\ox)(w),
\end{eqnarray*}
where $\ell\in\R_+$ is a Lipschitz constant of $\vt$ around $f(\ox)$ relative to its domain. This verifies the inequality ``$\ge$" in \eqref{dchain} and completes the proof of the theorem.
\end{proof}

Our next goal is to derive a subdifferential chain rule in the equality form for the limiting subdifferential in \eqref{sub} under the new metric subregularity qualification condition \eqref{mscq}. First we present the following lemma, which is of its own interest.\vspace*{-0.05in}

\begin{Lemma}[\bf extension of Lipschitz continuity]\label{extl} Let $\ph\colon\R^n\to\oR$ be a Lipschitz continuous function around $\ox$ relative to its domain with constant $\ell\in\R_+$. Then there exist a number $\ve>0$ and a function $\psi\colon\R^n\to\R$, which agrees with $\ph$ on $\B_\ve(\ox)\cap\dom\ph$ and which is Lipschitz continuous on $\R^n$ with the same constant $\ell$. If in addition $\ph$ is convex, then the function $\psi$ can be chosen to be convex as well.
\end{Lemma}\vspace*{-0.15in}
\begin{proof} The local Lipschitz continuity of $\ph$ relative to $\dom\ph$ gives us $\ve>0$ such that $\ph$ is Lipschitz continuous on the set $\O:=\B_\ve(\ox)\cap\dom\ph$ with constant $\ell$. Considering the function
$$
\psi(x):=\inf_{u\in\O}\big\{\ph(u)+\ell\|x-u\|\big\},\quad x\in\R^n,
$$
we can easily check (see, e.g., Rockafellar and Wets \cite[Exercise~9.12]{rw}) that $\psi$ agrees with $\ph$ on $\O$ while being Lipschitz continuous on $\R^n$ with the same constant $\ell$. Furthermore, the convexity of $\ph$ clearly yields the convexity of the function
$$
\theta(x,u):=\ph(u)+\delta_{\O}(u)+\ell\|x-u\|,\quad (x,u)\in\R^n\times\R^n,
$$
with respect to both variables. Having the representation $\psi(x)=\inf_{u\in\R^n}\theta(x,u)$, we deduce directly from the definition that $\psi$ is convex on $\R^n$.
\end{proof}\vspace*{-0.05in}

Now we are ready to establish the main result of this section providing the equality-type chain rule for limiting subgradients under the metric subregularity qualification condition \eqref{mscq}. Recall that a mapping $f\colon\R^n\to\R^m$ is {\em strictly differentiable} at $\ox$ with its strict derivative/Jacobian matrix $\nabla f(\ox)$ if we have
\begin{equation*}
\disp\lim_{x,u\to\ox}\frac{f(x)-f(u)-\nabla f(\ox)(x-u)}{\|x-u\|}=0.
\end{equation*}
{Observe that this differentiability notion lies between the usual Fr\'echet differentiability of $f$ at $\ox$ and its continuous differentiability around this point. Note also that the replacement in the next theorem of the conventional ${\cal C}^1$-smoothness of $f$ around $\ox$ by its merely strict differentiability at this point is necessary for applications to second-order analysis conducted in this paper. Indeed, our second-order analysis for \eqref{CS} is achieved under the twice differentiability of $f$ at the reference point. While implying the strict differentiability of $f$, this condition does not necessarily yield the ${\cal C}^1$-smoothness of this function.}\vspace*{-0.05in}

\begin{Theorem}[\bf equality chain rule for limiting subgradients under metric subregularity]\label{fchain} Let $f\colon\R^n\to\R^m$ be strictly differentiable at $\ox\in\R^n$, and let $\vt\colon\R^m\to\oR$ be convex, l.s.c.\ around $f(\ox)$, and Lipschitz continuous around this point relative to its domain with constant $\ell\in\R_+$. If MSQC \eqref{mscq} holds at $\ox$  with some constant $\kappa\in\R_+$, then the composition $\vt\circ f$ is lower regular at this point and we have the equality chain rule
\begin{equation}\label{1chain}
 \sub(\vt\circ f)(\ox)=\Hat\sub(\vt\circ f)(\ox)=\nabla f(\ox)^*\sub\vt\big(f(\ox)\big).
\end{equation}
\end{Theorem}\vspace*{-0.15in}
\begin{proof} Since our analysis is local around the points in question, there is no harm to suppose that $\vt$ is convex, l.s.c., and Lipschitz continuous relative to its entire domain. Remember that the inclusions
$$
\sub(\vt\circ f)(\ox)\supset\Hat\sub(\vt\circ f)(\ox)\supset\nabla f(\ox)^*\sub\vt\big(f(\ox)\big)
$$
always hold, we proceed with verifying the inclusion $\fsub(\vt\circ f)(\ox)\subset\nabla f(\ox)^*\sub\vt\big(f(\ox)\big)$, which justifies the claimed lower regularity of $\vt\circ f$ at $\ox$ together with the chain rule \eqref{1chain}. To this end, conclude from Lemma~\ref{extl} that there exists a Lipschitz continuous function $\psi\colon\R^m\to\R$ with constant $\ell$ such that $\vt=\psi+\delta_{\ss\dom\vt}$. This allows us to observe that
$$
\vt\circ f=\psi\circ f+\delta_{\ss\dom\vt}\circ f.
$$
It follows from the above representation of $\vt$ and the lower semicontinuity of $\vt$ that $\delta_{\ss\dom\vt}$ is l.s.c., and therefore $\dom\vt$ is closed. Using Lemma~2.1 from Gfrerer and Mordukhovich \cite{gm} gives us
$$
\sub(\delta_{\ss\dom\vt}\circ f)(\ox)\subset \nabla f(\ox)^*N_{\ss\dom\vt}\big(f(\ox)\big),
$$
which in fact holds as equality since $\dom\vt$ is convex. Employing the well-known chain and sum rules for Lipschitz continuous functions, we get
\begin{eqnarray}
\sub(\vt\circ f)(\ox)&\subset&\sub(\psi\circ f)(\ox)+\sub(\delta_{\ss\dom\vt}\circ f)(\ox)\nonumber\\
&\subset&\nabla f(\ox)^*\sub\psi\big(f(\ox)\big)+\nabla f(\ox)^*N_{\ss\dom\vt}\big(f(\ox)\big)\label{eq44}\\
&=&\nabla f(\ox)^*\sub\vt\big(f(\ox)\big)\nonumber,
\end{eqnarray}
where the first inclusion follows from \cite[Theorems~2.19]{m18} and the second one comes from \cite[Exercise~10.7]{rw}.
This verifies the claimed inclusion and thus completes the proof of the theorem.
\end{proof}\vspace*{-0.05in}

{All the assumptions on the outer function $\vt$ in Theorem~\ref{fchain} are automatically satisfied for fairly broad classes of extended-real-valued functions important in variational analysis and optimization; e.g., for the convex piecewise linear-quadratic function $\vt$ defined in \eqref{enlp}. It is worth mentioning that both composite functions considered in Example~\ref{fmr} enjoy MSQC \eqref{mscq}, and hence the chain rule \eqref{1chain} holds for them. As shown therein, the metric regularity qualification condition \eqref{bqc} fails for these composite functions.}\vspace*{0.05in}

The next corollary plays a significant role in deriving subsequent second-order results. It establishes the boundedness (with quantitative estimates) of dual elements under MSQC \eqref{mscq}. The latter is well known and rather easy to check under metric regularity.\vspace*{-0.05in}

\begin{Corollary}[\bf bounded multipliers]\label{subbound} Let $f\colon\R^n\to\R^m$ be strictly differentiable at $\ox\in\R^n$, and let $\vt\colon\R^m\to\oR$ be convex, l.s.c.\ around $f(\ox)$, and Lipschitz continuous around this point relative to its domain with constant $\ell\in\R_+$. If MSQC \eqref{mscq} holds at $\ox$  with some constant $\kappa\in\R_+$, then for every vector $v\in\sub(\vt\circ f)(\ox)$ there exists $\lambda\in\sub\vt(f(\ox))$ such that
\begin{equation}\label{subboundeq}
v=\nabla f(\ox)^*\lambda\;\mbox{ with }\;\|\lambda\|\le\ell+\kappa\|v\|+\kappa\ell\|\nabla f(\ox)\|.
\end{equation}
\end{Corollary}\vspace*{-0.15in}
\begin{proof} Assume without lost of generality that $\vt$ is Lipschitz continuous relative to its entire domain. Applying Lemma~\ref{extl} to the convex outer function $\vt$ in composition \eqref{CS}, we find a convex Lipschitz continuous function $\psi\colon\R^m\to\R$ such that $\vt=\psi+\delta_{\ss\dom\vt}$. Pick $v\in\sub(\vt\circ f)(\ox)$. It follows from \eqref{eq44} that there exist $\lambda_1\in\sub\psi(f(\ox))$ and $\lambda_2\in N_{\ss\dom\vt}(f(\ox))$ such that $v=\nabla f(\ox)^*(\lambda_1+\lambda_2)$. Since $\psi$ is Lipschitz continuous with the same constant $\ell$ due to Lemma~\ref{extl}, we have $\|\lambda_1\|\le\ell$. On the other hand, we can deduce from Lemma~2.1 in Gfrerer and Mordukhovich \cite{gm} that the following condition
$$
\|\lambda_2\|\le\kappa\|\nabla f(\ox)^*\lambda_2\|=\kappa\|v-\nabla f(\ox)^*\lambda_1\|.
$$
holds. Setting now $\lambda:=\lambda_1+\lambda_2$ leads us to
$$
\lambda\in\sub\psi\big(f(\ox)\big)+N_{\ss\dom\vt}\big(f(\ox)\big)=\sub\vt\big(f(\ox)\big).
$$
Furthermore, it follows from the above discussion that
\begin{eqnarray*}
\|\lambda\|=\|\lambda_1+\lambda_2\|\le\|\lambda_1\|+\kappa\|v-\nabla f(\ox)^*\lambda_1\|\le\ell+\kappa\|v\|+\kappa\ell\,\|\nabla f(\ox)\|,
\end{eqnarray*}
which readily verifies representation \eqref{subboundeq}.
\end{proof}\vspace*{-0.05in}

{We conclude this section by a consequence of the obtained chain rules, where the metric subregularity qualification condition \eqref{mscq} is automatically satisfied. These results cannot be deduced from the known qualification conditions formulated in Proposition~\ref{equi}(i,ii). While the chain rules below were mentioned in Rockafellar and Wets \cite[Exercise~10.22(b)]{rw},
the provided guide to the proof therein utilizes a number of rather involved results of variational analysis that are significantly different from our device.}\vspace*{-0.05in}

\begin{Corollary}[\bf chain rules for piecewise linear-quadratic functions]\label{chpwlq} Let $\ph\colon\R^n\to\oR$ be defined by $\ph(x):=\vt(Ax+a)$, where $\vt\colon\R^m\to\oR$ is a convex piecewise linear-quadratic function, $A$ is an $m\times n$ matrix, and $a\in\R^n$. Then for any point $x\in\dom\ph$ we have
\begin{equation*}
\d \ph(x)(w)=\d\vt(Ax+a)(Aw)\quad\mbox{and}\quad\sub\ph(x)=A^*\sub\vt(Ax+a).
\end{equation*}
\end{Corollary}\vspace*{-0.15in}
\begin{proof} Since $\vt$ is a convex piecewise linear-quadratic function, its domain is a polyhedral convex set. Then the Hoffman lemma tells us that MSQC \eqref{mscq} with $f(x):=Ax+a$ holds automatically at any point $x\in\dom\ph$. The claimed chain rules follow now from Theorems~\ref{regchain} and \ref{fchain}.
\end{proof}\vspace*{-0.2in}

\section{Fully Subamenable Functions}\label{subamen}
\sce

In this section we start the study and applications of major generalized differential constructions of second-order variational analysis. Our main attention is paid to the new class of fully subamenable functions. Prior to this, we recall yet another first-order subgradient notion and employ it to calculate the domain of second subderivatives \eqref{ssd} for a more general class of extended-real-valued functions.

Given $\ph\colon\R^n\to\oR$ and $\ox\in\dom\ph$, we say that $\ov\in\R^n$ is a {\em proximal subgradient} of $\ph$ at $\ox$ if there are positive numbers $\gamma$ and $r$ such that
\begin{equation}\label{prox9}
\ph(x)\ge\ph(\ox)+\langle\ov,x-\ox\rangle-\hbox{${r\over 2}$}\|x-\ox\|^2\quad\mbox{for all}\quad x\in\B_\gg(\ox).
\end{equation}
The set of all such $\ov$ is called the {\em proximal subdifferential} of $\ph$ at $\ox$ and is denoted by $\sub_p\ph(\ox)$.

The next theorem is important for its own sake and is also helpful for the subsequent derivations of second-order calculus rules and their applications.\vspace*{-0.05in}

\begin{Theorem}[\bf domain of second subderivatives]\label{properdom} Let $\ph\colon\R^n\to\oR$, and let $\ov\in\partial_p\ph(\ox)$ with $\ox\in\dom\ph$. The following assertions hold:

{\bf(i)} The second subderivative $\d^2\ph(\ox,\ov)$ is an l.s.c.\ function such that $\d^2\ph(\ox,\ov)(w)>-r\|w\|^2$ for some $r>0$ and all $w\in\R^n$. In particular, the function $\d^2\ph(\ox,\ov)$ is proper.

{\bf(ii)} If $T^2_{\ss{\epi\ph}}(\oz,q_w)\ne\emp$ with $\oz:=(\ox,\ph(\ox))$ and $q_w:=(w,\d\ph(\ox)(w))$ for all $w\in\R^n$ satisfying $\d\ph(\ox)(w)=\langle\ov,w\rangle$, then we have
\begin{equation}\label{domd2}
\dom\d^2\ph(\ox,\ov)=\big\{w\in\R^n\big|\;\d\ph(\ox)(w)=\langle\ov,w\rangle\big\}.
\end{equation}
\end{Theorem}\vspace*{-0.15in}
\begin{proof} The l.s.c.\ property of $\d^2\ph(\ox,\ov)$ was proved in Rockafellar and Wets \cite[Proposition~3.5]{rw}. To verify the lower estimate in (i), we take the triple $(\ov,r,\gg)$  from \eqref{prox9}, fix $w\in\R^n$, select $\ve>0$ with $\ve^2+\ve<\gg$, and pick any $u\in\B_\ve(w)$ and $t\in(0,\ve)$. If either $w=0$ or $\|w\|=1$, then $\ox+tu\in\B_\gg(\ox)$ and thus deduce from \eqref{prox9} that
\begin{equation*}
\Delta_t^2\ph(\ox,\ov)(u)=\frac{\ph(\ox+tu)-\ph(\ox)-t\langle\ov,u\rangle}{\sm t^2}\ge-r\|u\|^2.
\end{equation*}
This readily implies that for such $w$ we have the estimates
\begin{equation}\label{2dom}
\d^2\ph(\ox,\ov)(w)\ge-r\|w\|^2>-\infty.
\end{equation}
Since $w\mapsto\d^2\ph(\ox,\ov)(w)$ is positively homogeneous of degree $2$, we get
$$
\d^2\ph(\ox,\ov)(w)=\|w\|^2\d^2\ph(\ox,\ov)\Big(\frac{w}{\|w\|}\Big)>-r\|w\|^2\Big\|\frac{w}{\|w\|}\Big\|^2=-r\|w\|^2>-\infty,
$$
whenever $w\ne 0$, which verifies the estimates in \eqref{2dom} for all $w\in\R^n$. It is easy to observe that $\d^2\ph(\ox,\ov)(0)=0$, and thus $\d^2\ph(\ox,\ov)$ is a proper function.

To proceed next with the proof of (ii), note that the inclusion ``$\subset$" in \eqref{domd2} was already established in Rockafellar and Wets \cite[Proposition~13.5]{rw}. Let us derive the opposite inclusion for any fixed $w\in\R^n$ satisfying $\d\ph(\ox)(w)=\langle\ov,w\rangle$. By the assumed nonemptiness of the second-order set therein, we get
$(u,\alpha)\in T^2_{\scriptsize{\epi\ph}}(\oz,q_w)$ for the pair $(\oz,q_w)$ in the statement of the theorem. Then \eqref{2tan} gives us sequences $t_k\dn 0$ and $(u_k,\alpha_k)\to(u,\alpha)$ as $k\to\infty$ such that
\begin{equation*}
\big(\ox,\ph(\ox)\big)+t_k\big(w,\d\ph(\ox)(w)\big)+\sm t^2_{k}(u_k,\alpha_k)\in\epi\ph\;\mbox{ for all }\;k\in\N.
\end{equation*}
This tells us therefore that
\begin{equation*}
\frac{\ph\big(\ox+t_k w+\sm t^2_{k}u_k\big)-\ph(\ox)-t_k\d\ph(\ox)(w)}{\sm t^2_{k}}\le\alpha_k.
\end{equation*}
Denote $w_k:=w+\sm t_{k}u_k$ and deduce from $\d\ph(\ox)(w)=\langle\ov,w\rangle$ that
\begin{equation*}
\Delta_{t_k}^2\ph(\ox,\ov)(w_k)=\frac{\ph(\ox+t_k w_k) - \ph(\ox)-t_k\langle\ov,w_k\rangle}{\sm t^2_{k}}\le\alpha_k-\langle\ov,u_k\rangle.
\end{equation*}
Since $w_k\to w$ as $k\to\infty$, we arrive at the inequalities
\begin{equation*}
\d^2\ph(\ox,\ov)(w)\le\alpha-\langle\ov,u\rangle<\infty,
\end{equation*}
which show that $w\in\dom\d^2\ph(\ox,\ov)$ and hence complete the proof of the theorem.
\end{proof}

Next we introduce a fairly broad class of composite functions playing a crucial role in the rest of this paper. One of the nice properties of this class is that the sufficient condition for the validity of the domain formula \eqref{domd2} in Theorem~\ref{properdom} always holds for such functions, and thus we can use this formula for deriving the main second-order results.\vspace*{-0.05in}

\begin{Definition}[\bf fully subamenable functions]\label{full} We say that $\ph\colon\R^n\to\oR$ is {\sc fully subamenable} at $\ox\in\dom\ph$ if it there is a neighborhood $U$ of $\ox$ on which $\ph$ is represented as $\ph=\vt\circ f$, where $f\colon\R^n\to\R^m$ is twice differentiable at $\ox$, and where $\vt\colon\R^m\to\oR$ is convex piecewise linear-quadratic under the fulfillment of MSQC \eqref{mscq} at $\ox$.
\end{Definition}\vspace*{-0.05in}

The composition format of Definition~\ref{full} goes back to Rockafellar \cite{r88} who introduced the class of {\em fully amenable} functions, where MSQC \eqref{mscq} is replaced by the metric regularity qualification condition from Proposition~\ref{equi}(i), and where $f$ is assumed to be ${\cal C}^2$-smooth. Note that the latter metric regularity condition in this framework can be equivalently written in the Robinson constraint qualification form \eqref{gf07}. {We should point out here that both composite functions in Example~\ref{fmr} are clearly fully subamenable at the points in question. Moreover, the arguments in this example tell us that they are not fully amenable since the metric regularity qualification condition \eqref{bqc} fails for these composite functions.}\vspace*{0.05in}

Gfrerer and Mordukhovich \cite{gm} introduced the notions of (strongly, fully) {\em subamenable sets} of the type $\{x\in\R^n|\;f(x)\in\Th\}$ with the replacement of the metric regularity condition in Rockafellar \cite{r88} by the corresponding metric subregularity. This inspires the subamenability terminology in Definition~\ref{full}.\vspace*{0.03in}

{The result of Theorem~\ref{properdom}(ii) motivates the following notion of the critical cone for extended-real-value functions. Given $\psi\colon\R^n\to\oR$, define the {\em critical cone} of $\psi$ at $\ox$ for $\ov$ with $(\ox,\ov)\in\gph\sub\psi$ by
\begin{equation}\label{crit}
K_\psi(\ox,\ov):=\big\{w\in\R^n\big|\;\d\psi(\ox)(w)=\langle\ov,w\rangle\big\}.
\end{equation}
When $\psi=\dd_\O$ with $\O\subset\R^n$, we get $\d\psi(\ox)(w)=\dd_{T_\O(\ox)}(w)$ whenever $w\in\R^n$, which gives us the well-known critical cone construction for sets:
$$
K_\psi(\ox,\ov)=T_\O(\ox)\cap\{\ov\}^{\bot}.
$$
When $\vt\colon\R^m\to \oR$ a convex piecewise linear-quadratic function, it follows from \cite[Theorem~13.14]{rw} that the critical cone \eqref{crit} of $\vt$ at $(\ox,\ov)\in\gph\sub\vt$ can be equivalently described by
$$
K_\vt(\ox,\ov)=N_{\sub\vt(\ox)}(\ov).
$$}
Theorem~\ref{properdom}(ii) provides a useful sufficient condition under which the domain of the second subderivative and the critical cone agree. We show below that this condition is satisfied for fully subamenable functions. To proceed, let us first present the following chain rule for second-order tangent sets \eqref{2tan}. This result extends, with a different proof, the one in Rockafellar and Wets \cite[Proposition~13.13]{rw} from the fully amenable to fully subamenable setting.\vspace*{-0.05in}

\begin{Proposition}[\bf chain rule for second-order tangent sets]\label{tan-chain} Let $\ph\colon\R^n\to\oR$ be a fully subamenable composition at $\ox$. The following assertions hold:

{\bf(i)} The second-order tangent set $T^2_\O(\ox,w)$ is nonempty for any tangent vector $w\in T_\O(\ox)$, where the set $\O$ is defined in \eqref{mscq}.

{\bf(ii)} If $w\in T_\O(\ox)$ and the vector $z\in\R^n$ satisfy the inclusion
$$
\nabla f(\ox)z+\nabla^2 f(\bar x)(w,w)\in T^2_{\ss\dom\vt}\big(f(\bar x),\nabla f(\ox)w\big),
$$
then $z\in T^2_\O(\ox,w)$. Furthermore, there exist a number $\ve>0$ and an arc $\xi\colon[0,\ve]\to\O$ such that $\xi(0)=\ox$, $\xi_+'(0)=w$, and $\xi_+''(0)=z$.
\end{Proposition}\vspace*{-0.15in}
\begin{proof}
To verify (i), observe first that the inclusion $w\in T_\O(\ox)$ implies by MSQC \eqref{mscq} that $\nabla f(\ox)w\in T_{\ss\dom\vt}(f(\ox))$.
Since $\dom\vt$ is a polyhedral convex set, it follows that $f(\ox)+t\nabla f(\ox)w\in\dom\vt$ for all small $t>0$. Appealing now to \eqref{mscq}, we get for such $t$ that
\begin{eqnarray*}\label{chms}
{\rm dist}(\ox+t w;\O)&\le&\kappa\,{\rm dist}\big(f(\ox+t w);\dom\vt\big)\\
&\le&\kappa\,\big\|f(\ox+t w)-f(\ox)-t\nabla f(\ox) w\big\|=\frac{\kappa t^2}{2}\,\Big\|\nabla^2 f(\bar x)(w,w)+\frac{o(t^2)}{2 t^2}\Big\|.
\end{eqnarray*}
Thus there exists $u_t\in\O$ such that the parametric family of $z_t:=[u_t-\ox-t w]/\sm t^2$ is bounded for all small $t>0$.
Consequently, we find a sequence $t_k\dn 0$ for which $\ox+t_k w+\sm t_k^2 z_k\in\dom \vt $ and $z_k\to z$ as $k\to\infty$ with some $z\in\R^n$. By \eqref{2tan}, this yields $z\in T^2_\O(\ox,w)$ and hence verifies assertion (i).

Turning now to the proof of (ii), take $w\in T_\O(\ox)$ satisfying the relationships
$$
u:=\nabla f(\ox)z+\nabla^2 f(\bar x)(w,w)\in T^2_{\ss\dom\vt}\big(f(\bar x),\nabla f(\ox)w\big)=T_{T_{\ss\dom\vt}(f(\ox))}\big(\nabla f(\ox)w\big),
$$
where the equality is due to the polyhedrality of $\dom\vt$ and Proposition~13.12 in Rockafellar and Wets \cite{rw}. Hence we have
$f(\ox)+t\nabla f(\ox)w+\sm t^2 u\in\dom\vt$ for all small $t>0$. Then it follows from \eqref{mscq} that
\begin{eqnarray*}
{\rm dist}\big(\ox+t w+\sm t^2z;\O\big)&\le&\kappa\,{\rm dist}\big(f(\ox+t w+\sm t^2 z);\dom\vt\big)\\
&\le&\kappa\big\|f\big(\ox+t w+\sm t^2 z\big)-\big[f(\ox)+t\nabla f(\ox)w+\sm t^2 u\big]\big\|=o(t^2).
\end{eqnarray*}
Thus we find $z_t\in\O$ such that $\ox+t w+\sm t^2z-z_t=o(t^2)$  for all small $t>0$. Defining finally the arc $\xi(t):=z_t=\ox+t w+\sm t^2 z+o(t^2)$ verifies (ii) and completes the proof of the proposition.
\end{proof}\vspace*{-0.05in}

Now we are ready to derive the aforementioned result, which ensures the validity of the major assumption of Theorem~\ref{properdom} and of the critical cone formula for fully subamenable functions.\vspace*{-0.05in}

\begin{Theorem}[\bf critical cone for fully subamenable functions]\label{fsubam2} Let $\ph\colon\R^n\to\oR$ be a fully subamenable function at $\ox$ in the notation of Theorem~{\rm\ref{properdom}}. Then we have $\sub\ph(\ox)=\sub_p\ph(\ox)$ and
$T^{2}_{\scriptsize{\epi\ph}}(\oz,q_w)\ne\emp$ for any $w\in T_\O(\ox)$. Consequently
\begin{equation}\label{domdf}
K_\ph(\ox,\ov)=\dom\d^2\ph(\ox,\ov)\;\mbox{ whenever }\;\ov\in\sub\ph(\ox).
\end{equation}
\end{Theorem}\vspace*{-0.15in}
\begin{proof} To verify the first statement, observe that we always have $\sub_p\ph(\ox)\subset\sub\ph(\ox)$ and proceed with the proof of the opposite inclusion. Pick $\ov\in\sub\ph(\ox)$ and find by \eqref{1chain} a vector $\olambda\in\sub\vt(f(\ox))$ with $ \ov=\nabla f(\ox)^*\olambda$. It follows from the twice differentiability of $f$ at $\ox$ that
\begin{equation*}
f(x)-f(\ox)=\nabla f(\ox)(x-\ox)+\nabla^2 f(\ox)(x-\ox,x-\ox)+o(\|x-\ox\|^2).
\end{equation*}
Combining it with the convexity of $\vt$ gives us $\ve>0$ such that
\begin{eqnarray*}
\vt\big(f(x)\big)-\vt\big(f(\ox)\big)\ge\langle\olambda,f(x)-f(\ox)\rangle&\ge&\langle\olambda,\nabla f(\ox)(x -\ox)\rangle-\|\olambda\|(1+\|\nabla^2 f(\ox)\|)\|x-\ox\|^2\\
&=&\langle\ov,x-\ox\rangle-\|\olambda\|(1+\|\nabla^2 f(\ox)\|)\|x-\ox\|^2
\end{eqnarray*}
for all $x\in\B_{\ve}(\ox)$. This shows that $\ov\in\sub_p\ph(\ox)$ and thus justifies the claimed inclusion.

To verify next that $T^{2}_{\scriptsize{\epi\ph}}(\oz,q_w)\ne\emp$, fix $w\in T_\O(\ox)$ and get from Proposition~\ref{tan-chain}(i) that $T^2_\O(\ox,w)\ne\emp$. Picking $u\in T^2_\O(\ox,w)$, we find a sequence $t_k\dn 0$ such that
\begin{equation*}\label{fs1}
\ox+t_k w+\sm t_k^2u+o(t_k^2)\in\O\;\mbox{ for all }\;k\in\N.
\end{equation*}
Denote $w_k:= w+\frac{1}{2} t_k u+\frac{o(t_k^2)}{t_k}$ and deduce from the definition of $\O$ in \eqref{mscq} that
$f(\ox+t_kw_k)\in\dom\vt$, which gives us the relationships
\begin{equation*}\label{fs2}
f(\ox+t_k w_k)=f(\ox)+t_k\nabla f(\ox)w_k+\sm t_k^2\nabla^2 f(\ox)(w_k,w_k)+o(t_k^2)\in\dom\vt.
\end{equation*}
Since $\vt$ is piecewise linear quadratic, we have $\dom\vt=\cup_{j=1}^{s}\O_j$, where each $\O_j$ is a polyhedral convex set. Passing to a subsequence if necessary, choose an index $i\in\{1,\ldots,s\}$ for which
\begin{equation*}
y_k:=f(\ox+t_k w_k)=f(\ox)+t_k\nabla f(\ox)w_k+\sm t_{k}^{2}\nabla^2 f(\ox)(w_k,w_k)+o(t_{k}^{2})\in\O_{i}
\end{equation*}
whenever $k\in\N$. Remembering that $\vt(y)=\alpha_i+\langle a_i,y\rangle+\sm\langle A_i y,y\rangle$ for all $y\in\O_i$ yields
\begin{eqnarray*}
\ph(\ox+t_k w_k)&=&\vt\big(f(\ox+t_k w_k)\big)=\vt(y_k)=\alpha_i+\langle a_i,y_k\rangle+\sm\langle A_iy_k,y_k\rangle\\
&=&\alpha_i+\langle a_i,f(\ox)\rangle+\sm\langle A_i f(\ox),f(\ox)\rangle+t_k\langle A_if(\ox)+a_i,\nabla f(\ox)w_k\rangle\\
&&+\sm t_{k}^2\big[\langle A_i\nabla f(\ox)w_k,\nabla f(\ox)w_k\rangle+\langle a_i,\nabla^2 f(\ox)(w_k,w_k)\rangle\\
&&+2\langle A_if(\ox),\nabla^2f(\ox)(w_k,w_k)\rangle\big]+o(t_{k}^{2})\\
&=&\ph(\ox)+t_k\langle A_if(\ox)+a_i,\nabla f(\ox)w\rangle+\sm\langle A_if(\ox)+a_i,\nabla f(\ox)u\rangle+o(t_{k}^{2})\\
&&+\sm t_{k}^2\big[\langle A_i\nabla f(\ox)w_k,\nabla f(\ox)w_k\rangle+\langle a_i,\nabla^2 f(\ox)(w_k,w_k)\rangle\big],
\end{eqnarray*}
where in the last equality we used that $w_k=w+\sm t_k u+\frac{o(t_{k}^{2})}{t_k}$. This tells us that
\begin{eqnarray*}
\big(\ox+t_k w_k,\ph(\ox+t_k w_k)\big)&=&\big(\ox,f(\ox)\big)+t_k\Big(w,\big\langle A_if(\ox)+a_i,\nabla f(\ox)w\big\rangle\Big)\\
&&+\sm t_{k}^2\Big(u,\big\langle A_i\nabla f(\ox)w_k,\nabla f(\ox)w_k\big\rangle+\big\langle a_i,\nabla^2 f(\ox)(w_k,w_k)+\nabla f(\ox)u\big\rangle\\
&&+2\big\langle A_if(\ox),\nabla^2 f(\ox)(w_k,w_k)+\nabla f(\ox)u\big\rangle\Big)+o(t_{k}^{2})\in\gph \ph\subset\epi\ph.
\end{eqnarray*}
It follows from the chain rule for subderivatives in Theorem~\ref{regchain} and from the subderivative representation for piecewise linear-quadratic functions in \eqref{dfs2} that
$$
\d\ph(\ox)(w)=\d\vt\big(f(\ox)\big)\big(\nabla f(\ox)w\big)=\big\langle A_if(\ox)+a_i,\nabla f(\ox)w\big\rangle.
$$
Combining these and letting $k\to\infty $, we arrive at $(u,p)\in T^{2}_{\scriptsize{\epi\ph}}(\ox,q_w)$ with
\begin{equation*}
p:=\big\langle A_i\nabla f(\ox) w,\nabla f(\ox)w\big\rangle+\big\langle a_i,\nabla^2 f(\ox)(w,w)+\nabla f(\ox)u\big\rangle+2\langle A_if(\ox),\nabla^2f(\ox)(w,w)+\nabla f(\ox)u\rangle,
\end{equation*}
which justifies the claimed nonemptiness of the second-order tangent set.

It remains to verify the critical cone formula \eqref{domdf}. Pick $\ov\in\sub\ph(\ox)$ and fix a vector $w\in\R^n$ satisfying $\d\ph(\ox)(w)=\langle\ov,w\rangle$. Furthermore, it follows from the assumed MSQC \eqref{mscq} and from the given formula for the domain of $\d\vt(f(\ox))$ in  \eqref{dfs} that
$$
T_\O(\ox)=\big\{w\in\R^n\big|\;\nabla f(\ox)w\in T_{\ss\dom\vt}\big(f(\ox)\big)\big\}=\big\{w\in\R^n\big|\;\nabla f(\ox)w\in\dom\d\vt\big(f(\ox)\big)\big\},
$$
which yields $w\in T_\O(\ox)$, and hence $T^{2}_{\scriptsize{\epi\ph}}(\oz,q_w)\ne\emp$ as shown above. Appealing now to Theorem~\ref{properdom}(ii) justifies \eqref{domdf} and thus completes the proof.
\end{proof}\vspace*{-0.05in}

Theorem~\ref{fsubam2} plays a key role in the {\em variational approach} to second-order calculus and applications of fully subamenable compositions developed in the subsequent sections. This is largely due to the next theorem that establishes the existence of optimal solutions to a special class of constrained optimization problems constructed via the second-order data of fully subamenable functions. To define such a problem for a given fully subamenable function $\ph$ at $\ox$, deduce first from the subdifferential chain rule \eqref{1chain} and Proposition~\ref{nes} that $\sub\ph(\ox)\ne\emp$. Pick $\ov\in\sub\ph(\ox)$ and define the {\em multiplier set} associated with $(\ox,\ov)$ by
\begin{equation}\label{lagn}
\Lambda(\ox,\ov):=\big\{\lambda\in\R^m\big|\;\nabla f(x)^*\lambda=\ov,\;\lambda\in\partial\vt\big(f(\ox)\big)\big\}.
\end{equation}
The imposed MSQC \eqref{mscq} ensures that this set is nonempty. Moreover, $\Lambda(\ox,\ov)$ is a polyhedral convex set since the function $\vt$ is convex and piecewise linear-quadratic. Fix now a vector $w\in\R^n$ and consider the constrained {\em optimization problem}:
\begin{equation}\label{du}
\max_{\lambda\in\R^{m}}\;\;\left\langle\lambda,\nabla^2 f(\ox)(w,w)\right\rangle+\d^2\vt\big(f(\ox),\lm\big)\big(\nabla f(\ox)w\big)\quad\mbox{subject to}\quad \lambda\in\Lambda(\ox,\ov).
\end{equation}
We need to use the following result on second subderivatives of convex piecewise linear-quadratic functions that can be extracted from the proof of Proposition~13.9 in Rockafellar and Wets \cite{rw}: Let $\vt\colon\R^m\to\oR$ be given in form \eqref{PWLQ} with $\oy\in\dom\vt$. Then $\vt$ is properly twice epi-differentiable at $\oy$ for any $\ou\in\sub\vt(\oy)$ with the representations
\begin{equation}\label{pwfor}
\dom\d^2\vt(\oy,\ou)=\bigcup_{i\in I(\oy)}T_{\O_i}(\oy)\cap\{\ou_i\}^{\bot},\;\d^2\vt(\oy,\ou)(w)
=\begin{cases}
\la A_i w,w\ra&\mbox{if}\;\;w\in T_{\O_i}(\oy)\cap\{\ou_i\}^{\bot},\\
\infty&\mbox{otherwise},
\end{cases}
\end{equation}
where $\ou_i:=\ou-A_i\oy-a_i$, and where the index set $I(\oy)$ is taken from \eqref{dfs}. Now we are ready to establish the aforementioned existence result.\vspace*{-0.05in}

\begin{Theorem}[\bf existence of optimal solutions along critical directions]\label{duno} Let $\ph\colon\R^n\to\oR$ be a fully subamenable function at $\ox$, and let $\ov\in\sub\ph(\ox)$. Then for any critical direction $w\in K_\ph(\ox,\ov)$ we have:

{\bf(i)} There exists an optimal solution to problem \eqref{du}.

{\bf(ii)} Denoting ${\cal A}:=\big\{\lm\in\sub\vt(f(\ox))\big|\;\d\vt(f(\ox))(\nabla f(\ox)w)=\la\nabla f(\ox)w,\lm\ra\big\}$, there exists an optimal solution to the modified problem
\begin{equation}\label{du2}
\max_{\lambda\in\R^{m}}\;\left\langle\lambda,\nabla^2f(\ox)(w,w)\right\rangle+\d^2\vt\big(f(\ox),\lm\big)\big(\nabla f(\ox)w\big)\;\mbox{ subject to }\;\nabla f(x)^*\lambda=\ov,\;\lambda\in{\cal A}.
\end{equation}
Moreover, the sets of optimal solutions to problems \eqref{du} and \eqref{du2} coincide.
\end{Theorem}\vspace*{-0.15in}
\begin{proof} Fix $w\in\R^n$ and pick any $\lm\in\Lambda(\ox,\ov)$. Arguing similarly to the beginning of the proof of Theorem~13.14 in Rockafellar and Wets \cite{rw} (this part of the proof can be carried out for fully subamenable functions, not just for fully amenable ones as assumed therein), we get
\begin{equation}\label{ine2}
\d^2\ph(\ox,\ov)(w)\ge\d^2\vt\big(f(\ox),\lambda\big)\big(\nabla f(\ox)w\big)+\left\langle\lambda,\nabla^2f(\ox)(w,w)\right\rangle.
\end{equation}
Take now $w\in K_\ph(\ox,\ov)$, which is equivalent by MSQC \eqref{mscq} to $\nabla f(\ox)w\in K_\vt(f(\ox),\lm)$. It follows from \eqref{domdf} that both numbers $\d^2\ph(\ox,\ov)(w)$ and $\d^2\vt(f(\ox),\lambda)(\nabla f(\ox) w)$ are finite, and hence
$$
\sup_{\lambda\in\Lambda(\ox,\ov)}\big\{\left\langle\lambda,\nabla^2f(\ox)(w,w)\right\rangle+\d^2\vt\big(f(\ox),\lm\big)(\nabla f(\ox)w)\big\}\le\d^2\ph(\ox,\ov)(w)<\infty.
$$
This ensures that the optimal value of \eqref{du} is finite. Furthermore, we observe from \eqref{pwfor} that for any $w\in K_\ph(\ox,\ov)$ the second subderivative $\d^2\vt(f(\ox),\lm)(\nabla f(\ox)w)$ is actually independent of $\lm$. This tells us that problem \eqref{du} is a {\em linear program}, where the optimal value is finite. This yields the existence of optimal solutions to \eqref{du}, which verifies (i).

To check (ii), we deduce from $w\in K_\ph(\ox,\ov)$ and formula \eqref{domdf} in Theorem~\ref{fsubam2} that $\d\ph(\ox)(w)=\langle\ov,w\rangle$. Employing the chain rule for subderivatives from \eqref{dchain} gives us
$$
\d\vt\big(f(\ox)\big)\big(\nabla f(\ox)w\big)=\la\nabla f(\ox)w,\lm\ra\quad\mbox{for all}\quad\lm\in\Lambda(\ox,\ov).
$$
This means that the constraint $\d\vt(f(\ox))(\nabla f(\ox)w)=\la\nabla f(\ox)w,\lm\ra$ is an {\em implicit} constraint for problem \eqref{du}, and therefore the sets of feasible solutions to problems \eqref{du} and \eqref{du2} agree. This verifies by using (i) that for any $w\in K_\ph(\ox,\ov)$ the sets of optimal solutions to problems \eqref{du2} and \eqref{du} are the same.
\end{proof}\vspace*{-0.05in}

Besides the optimization problem \eqref{du} and its equivalent subderivative form \eqref{du2}, in the next section we deal with its {\em dual} problem constructed by using parabolic subderivatives. Given $\psi\colon\R^n\to\oR$ finite at $\ox$ and given $w\in\R^n$ where $\d\psi(\ox)(w)$ is finite, the {\em parabolic subderivative} of $\psi$ at $\ox$ for $w$ with respect to $z\in\R^m$ is defined by
\begin{equation}\label{par-subder}
\d^2\psi(\bar x)(w;z):=\liminf_{\substack{t\dn 0\\u\to z}}\dfrac{\psi(\ox+tw+\sm t^2 u)-\psi(\ox)-t\d\psi(\ox)(w)}{\sm t^2}.
\end{equation}

Let us first summarize some well-known properties of parabolic subderivatives for the class of convex piecewise linear-quadratic functions.\vspace*{-0.05in}

\begin{Proposition}[\bf properties of parabolic subderivatives and Fenchel conjugates]\label{parp} Let $\psi\colon\R^n\to\oR$ be convex piecewise linear-quadratic with $\ox\in\dom\psi$, and let $w\in\dom\d\psi(\ox)$. Then the function $z\mapsto\d^2\psi(\bar x)(w;z)$ is proper, l.s.c., and convex piecewise linear with its Fenchel conjugate calculated by
\begin{equation}\label{conpi}
v\mapsto\begin{cases}
-\d^2\psi(\ox,v)(w)&\mbox{if}\;v\in\sub\psi(\ox)\;\mbox{ with }\;\d\psi(\ox)(w)=\la w,v\ra,\\
\infty&\mbox{otherwise}.
\end{cases}
\end{equation}
Furthermore, we have the equivalence
\begin{equation}\label{dompa}
\d^2\psi(\bar x)(w;z)<\infty\iff z\in T_{T_{\ss\dom\psi}(\oy)}(w)=T^2_{\ss\dom\psi}(\ox,w).
\end{equation}
\end{Proposition}\vspace*{-0.1in}
\begin{proof} It is easy to derive the listed general properties from the definitions of parabolic subderivatives and piecewise linear-quadratic functions. The Fenchel conjugate formula \eqref{conpi} is due to Rockafellar \cite[Proposition~3.5]{r88}, while the equivalence \eqref{dompa} is discussed in Rockafellar and Wets \cite[Exercise~13.61]{rw}.
\end{proof}\vspace*{-0.05in}

Employing finally Theorem~\ref{duno} together with Proposition~\ref{parp} leads us to duality relationships for piecewise linear-quadratic programs that extend the classical ones in linear programming by taking into account specific structures of problems \eqref{du} and \eqref{du2}.\vspace*{-0.05in}

\begin{Corollary}[\bf duality relationships along critical directions]\label{pidu} Let $\ph\colon\R^n\to\oR$ be a fully subamenable function at $\ox$, and let $\ov\in\sub\ph(\ox)$. Then for any $w\in K_\ph(\ox,\ov)$ we have:

{\bf(i)} The dual problem of the piecewise linear-quadratic program \eqref{du2} is given by
\begin{equation}\label{pri}
\min_{z\in\R^{n}}\;-\left\langle\ov,z\right\rangle+\d^2\vt\big(f(\ox)\big)\big(\nabla f(\ox)w;\nabla f(\ox)z+\nabla^2f(\ox)(w,w)\big),
\end{equation}
and it admits an optimal solution.

{\bf(ii)} The optimal value of \eqref{du2} coincides with the optimal value of \eqref{pri}.
\end{Corollary}\vspace*{-0.15in}
\begin{proof} The duality in (i) follows from Proposition~\ref{parp}, while the existence of optimal solutions to \eqref{pri} and the claim in (ii) are consequences of Theorem~\ref{duno}.
\end{proof}\vspace*{-0.2in}

\section{Second Subderivatives of Fully Subamenable Functions}\label{2subamen}
\sce\vspace*{-0.1in}

The main goal of this section is to establish twice epi-differentiability of every fully subamenable function. For fully amenable functions it was done by Rockafellar in \cite{r88} with the detailed proof given in Theorem~13.14 of his book with Wets \cite{rw}. Our approach here is significantly different and in fact much simpler even for fully amenable functions; see the comments after the proof. One of the new ingredients is involving into the proof the parabolic subderivatives \eqref{par-subder}. Furthermore, we obtain precise formulas for second subderivatives of fully subamenable functions that are expressed entirely via the given data.

The following lemma of its own interest is useful in the proof of the main result of this section. It reveals an important second-order property of outer functions that appear in fully subamenable compositions.\vspace*{-0.05in}

\begin{Lemma}[\bf parabolic subderivatives of piecewise linear-quadratic functions]\label{parsub} Let $\vt\colon\R^m\to\oR$ be a convex piecewise linear-quadratic function with representation \eqref{PWLQ}. Then for any $\oy\in\dom\vt$, $w\in\dom\d\vt(\oy)$, and $z\in T^2_{\ss\dom\vt}(\oy,w)$ we have
\begin{equation}\label{par-pwlq}
\d^2\vt(\bar y)(w,z)=\lim_{t\dn 0}\dfrac{\vt(\oy+tw+\sm t^2z)-\vt(\oy)-t\d\vt(\oy)(w)}{\sm t^2}.
\end{equation}
\end{Lemma}\vspace*{-0.15in}
\begin{proof} To verify the claimed representation \eqref{par-pwlq}, remember that $z\in T^2_{\ss\dom\vt}(\oy,w)$, and thus it follows from the equivalence in \eqref{dompa} that $\d^2\vt(\bar y)(w,z)$ is finite. Thus by the definition of the parabolic subderivative, we find sequences $t_k\dn 0$ and $z_k\to z$ such that
$$
\lim_{k\to\infty}\dfrac{\vt\big(\oy+t_kw+\sm t_k^2z_k\big)-\vt(\oy)-t_k\d\vt(\oy)(w)}{\sm t_k^2}=\d^2\vt(\bar y)(w,z).
$$
Since $\d^2\vt(\bar y)(w,z)$ is finite, we have $\oy+t_kw+\sm t_k^2 z_k\in\dom\vt=\cup_{i=1}^s\O_i$ for all $k$ sufficiently large. Passing to a subsequence if necessary, find an index $i_0\in\{1,\ldots,s\}$ with $\oy+t_kw+\sm t_k^2 z_k\in\O_{i_0}$ for such large $k\in\N$. It ensures therefore that $i_0\in J(\oy,w)$.
Appealing now to \eqref{PWLQ} gives us
\begin{eqnarray*}
\vt\big(\oy+t_k(w+\sm t_k z_k)\big)-\vt(\oy)&=&t_k\big\la A_i\oy+a_i,w+\sm t_k z_k\big\ra+\sm t_k^2\big\la A_i\big(w+\sm t_k z_k\big),\big(w+\sm t_k z_k\big)\big\ra\\
&=& t_k\d\vt(\oy)(w)+\sm t_k^2\la A_i\oy+a_i,z_k\ra+\sm t_k^2\la w,A_iw\ra+o(t_k^2).
\end{eqnarray*}
It follows from $i_0\in J(\oy,w)$ that $i_0\in I(\oy)$ and thus $\d\vt(\oy)(w)=\la A_{i_0}\oy+a_{i_0},w\ra$ by \eqref{dfs2}.
Combining these in turn brings us to the equality
\begin{equation}\label{par3}
\lim_{k\to\infty}\dfrac{\vt\big(\oy+t_kw+\sm t_k^2 z_k\big)-\vt(\oy)-t_k\d\vt(\oy)(w)}{\sm t_k^2}=\la A_{i_0}\oy+a_{i_0},z\ra+\la w,A_{i_0}w\ra.
\end{equation}
The above arguments verify  that $z\in T^2_{\O_{i_0}}(\oy,w)$. The polyhedrality of $\O_{i_0}$ ensures the existence of $\ve>0$ with $\oy+tw+\sm t^2z\in\O_{i_0}$ for all $t\in[0,\ve]$. Arguing as in the proof of \eqref{par3} yields
\begin{equation*}
\lim_{t\dn 0}\dfrac{\vt\big(\oy+tw+\sm t^2 z\big)-\vt(\oy)-t\d\vt(\oy)(w)}{\sm t^2}=\la A_{i_0}\oy+a_{i_0},z\ra+\la w,A_{i_0}w\ra.
\end{equation*}
This along with \eqref{par3} confirms that both sides of the equality in \eqref{par-pwlq} agree with each other for any second-order tangent vector $z\in T^2_{\ss\dom\vt}(\oy,w)$.
\end{proof}\vspace*{-0.05in}

Now we are ready to prove the twice epi-differentiability of fully subamenable functions and derive an explicit formula for their second subderivatives.\vspace*{-0.05in}

\begin{Theorem}[\bf twice epi-differentiability of fully subamenable functions]\label{sochain} Let $\ph\colon\R^n\to\oR$ be a fully subamenable function at $\ox\in\dom\ph$. Then $\ph$ is properly twice epi-differentiable at $\ox$ for every $\ov\in\sub\ph(\ox)$, and its second subderivative \eqref{ssd} is calculated by
\begin{equation}\label{tedfu}
\d^2\ph(\ox,\ov)(w)=\max_{\lambda\in\Lambda(\ox,\ov)}\big\{\d^2\vt\big(f(\ox),\lambda\big)\big(\nabla f(\ox)w\big)+\left\langle\lambda,\nabla^2f(\ox)(w,w)\right\rangle\}
\end{equation}
for all $w\in\R^n$, where the set of multiplies $\Lambda(\ox,\ov)$ is taken from \eqref{lagn}.
\end{Theorem}\vspace*{-0.15in}
\begin{proof} The inequality ``$\ge$" in \eqref{tedfu} is given in \eqref{ine2}. Now we proceed with the simultaneous verification of the opposite inequality in \eqref{tedfu} and the twice epi-differentiability of $f$ in \eqref{df02}. It only suffices to prove these relationships for critical directions $w\in K_\ph(\ox,\ov)$. Indeed, for $w\notin K_\ph(\ox,\ov)$, which is equivalent to $\nabla f(\ox)w\notin K_\vt(f(\ox),\lm)$ whenever $\lm\in\Lambda(\ox,\ov)$, both sides of \eqref{tedfu} become $\infty$. To obtain \eqref{df02} for this case, pick a sequence $t_k\dn 0$ and then let the sequence $w_k:=w$ for any $k$. It is easy to observe that \eqref{df02} holds for the aforementioned sequence. Fix further $w\in K_\ph(\ox,\ov)$ and pick any $\ov\in\sub\ph(\ox)$. This together with the chain rule for subderivatives \eqref{dchain} ensures the equalities
\begin{equation}\label{jg02}
\d\ph(\ox)(w)=\d\vt\big(f(\ox)\big)\big(\nabla f(\ox)w\big)=\la\ov,w\ra.
\end{equation}
Corollary~\ref{pidu}(i) tells us that the piecewise linear-quadratic program \eqref{pri} admits an optimal solution denoted by $\oz$, and its optimal value is finite. Hence the parabolic subderivative $\d^2\vt(f(\ox))\big(\nabla f(\ox)w;\nabla f(\ox)\oz+\nabla^2f(\ox)(w,w)\big)$ is also finite. This together with
\eqref{dompa} implies that
\begin{equation}\label{so2}
\ou:=\nabla^2f(\ox)(w,w)+\nabla f(\ox)\oz\in T_{\ss\dom\vt}^2\big(f(\ox),\nabla f(\ox)w\big).
\end{equation}
Since $\dom \vt$ is a polyhedral convex set, we find a number $\dd>0$ such that
\begin{equation}\label{jg01}
f(\ox)+t\nabla f(\ox)w+\sm t^2\ou\in\dom\vt\quad\mbox{for all}\quad t\in[0,\dd].
\end{equation}
Moreover, we conclude from \eqref{so2} and Proposition~\ref{tan-chain}(ii) that $\oz\in T_\O^2(\ox,w)$. The latter proposition also ensures the existence of a number $\ve\in(0,\dd)$ and an arc $\xi\colon[0,\ve]\to\O$ for which we have
\begin{equation}\label{df03}
\xi(0)=\ox,\;\xi'_{+}(0)=w,\;\mbox{ and }\;\xi''_{+}(0)=z.
\end{equation}
Define now $w_t:=\disp\frac{\xi(t)-\xi(0)}{t}$ for all $t\in[0,\ve]$ and observe that $\ox+tw_t=\xi(t)\in\O$ for such $t$. It follows from the second equality in \eqref{df03} that $w_t\to w$ as $t\dn 0$. Thus for all $t\in[0,\ve]$ we deduce from the relationships in \eqref{jg02}--\eqref{df03} that
\begin{eqnarray}\label{seq5}
\begin{array}{ll}
&\Delta_t^2\ph(\ox,\ov)(w_t)=\dfrac{\ph(\ox+tw_t)-\ph(\ox)-t\langle\ov,w_t\rangle}{\sm t^2}=\dfrac{\vt\big(f(\ox+tw_t)\big)-\vt\big(f(\ox)\big)-t\langle\ov,w_t\rangle}{\sm t^2}\\
&=\dfrac{\vt\big(f(\xi(t))\big)-\vt\big(f(\ox)\big)-t\,\d\vt\big(f(\ox)\big)\big(\nabla f(\ox)w\big)}{\sm t^2}-\Big\la\ov,\dfrac{\xi(t)-\xi(0)-t w}{\sm t^2}\Big\ra\\
&=\dfrac{\vt\big(f(\ox)+t\nabla f(\ox)w+\sm t^2\ou\big)-\vt\big(f(\ox)\big)-t\,\d\vt\big(f(\ox)\big)\big(\nabla f(\ox)w\big)}{\sm t^2}\\
&+\dfrac{\vt\big(f(\xi(t))\big)-\vt\big(f(\ox)+t\nabla f(\ox)w+\sm t^2\ou\big)}{\sm t^2}-\Big\la\ov,\dfrac{\xi(t)-\xi(0)-t\xi'_{+}(0)}{\sm t^2}\Big\ra.
\end{array}
\end{eqnarray}
Looking at the last equality in \eqref{seq5}, we see that the first term therein converges to $\d^2\vt(f(\ox))(\nabla f(\ox)w;\ou)$ as $t\dn 0$ due to Lemma~\ref{parsub}. The third term clearly converges to $-\la\ov,\oz\ra$. Turing to the second term in this equality, remember that $\vt$ is Lipschitz continuous relative to its domain, which implies by \eqref{jg01} and $f(\xi(t))\in\dom\vt$ that the second term converges to zero since
$$
\frac{f\big(\xi(t)\big)-f(\ox)-t\nabla f(\ox)w}{\sm t^2}\to\ou\quad\mbox{as}\quad t\dn 0.
$$
Getting all the above together, we arrive at the equalities
\begin{eqnarray*}
\lim_{t\dn 0}\Delta_t^2\ph(\ox,\ov)(w_t)&=&\d^2\vt\big(f(\ox)\big)\big(\nabla f(\ox)w;\ou\big)-\la\ov,\oz\ra\\
&=&\d^2\vt\big(f(\ox)\big)\big(\nabla f(\ox)w;\nabla f(\ox)\oz+\nabla^2f(\ox)(w,w)\big)-\la\ov,\oz\ra\\
&=&\max_{\lambda\in\Lambda(\ox,\ov)}\big\{\d^2\vt\big(f(\ox),\lambda\big)\big(\nabla f(\ox)w\big)+\left\langle\lambda,\nabla^2f(\ox)(w,w)\right\rangle\big\},
\end{eqnarray*}
where the last one comes from Proposition~\ref{pidu}(ii). This verifies the inequality ``$\le$"  in \eqref{tedfu} as well as the convergence in \eqref{df02}, and thus completes the proof of the theorem.
\end{proof}\vspace*{-0.05in}

As mentioned above, the results of Theorem~\ref{sochain} extend those in Rockafellar and Wets \cite[Theorem~13.14]{rw} obtained under the metric regularity qualification condition in form \eqref{gf07}. Our proof, based on metric subregularity in \eqref{mscq}, is largely different and essentially simpler than the one from \cite[Theorem~13.14]{rw}. The major difference is that we use in the proof another pair of primal-dual problems and involve parabolic subderivatives. Our approach allows us to deal with more general frameworks and applications, which is the main subject of our subsequent research \cite{mms}.\vspace*{0.05in}

Remember that the equivalent optimization problems \eqref{du} and \eqref{du2} are in fact problems of {\em linear programming} due to the second subderivative calculation \eqref{pwfor}. Next we derive equivalent representations for the second subderivative of $\ph$ in \eqref{tedfu}, which are useful for proto-derivative calculations in Section~\ref{gra-der}.\vspace*{-0.05in}

\begin{Corollary}[\bf chain rules for second subderivatives of fully subamenable compositions]\label{cso0} In the framework and under the assumptions of Theorem~{\rm\ref{sochain}}, take any $\bar\lm\in\Lambda(\ox,\ov)$. Then we have the second-order chain rules:

{\bf(i)} The second subderivative of $\ph$ at $\ox$ for $\ov$ is represented by
\begin{equation}\label{cso}
\d^2\ph(\ox,\ov)(w)=\d^2\vt\big(f(\ox),\olambda\big)\big(\nabla f(\ox)w\big)+\max_{\lambda\in\Lambda(\ox,\ov)}\big\langle\lambda,\nabla^2f(\ox)(w ,w)\big\rangle.
\end{equation}

{\bf(ii)} There exits $\bar r>0$ such that for any $r>\bar r$ and any $w\in\R^n$ we have
\begin{equation}\label{cso2}
\d^2\ph(\ox,\ov)(w)=\d^2\vt\big(f(\ox),\olambda\big)\big(\nabla f(\ox)w\big)+\max_{\lambda\in[\Lambda(\ox,\ov)\,\cap\,r\B]}\big\langle\lambda,\nabla^2f(\ox)(w,w) \big\rangle.
\end{equation}
\end{Corollary}\vspace*{-0.1in}
\begin{proof} We get from the critical cone definition and the subderivative chain rule under MSQC \eqref{mscq} that
\begin{equation}\label{pj98}
w \in K_\ph(\ox,\ov)\iff\nabla f(\ox)w\in K_\vt\big(f(\ox),\lm\big)\;\mbox{ for all }\;\lm\in\Lambda(\ox,\ov).
\end{equation}
If $w\notin K_\ph(\ox,\ov)$, then both sides in \eqref{cso}  become $\infty$, and so the equality holds therein. Similarly we get \eqref{cso2} for any $r>0$ in this case.

Consider now the case where $w\in K_\ph(\ox,\ov)$. To verify (i), deduce from \eqref{pwfor} that the second subderivative $\d^2\vt(f(\ox),\olambda)(\nabla f(\ox)w)$ is actually {\em independent} of $\bar\lm$ for such critical directions $w$, and so the sets on the right-hand sides in \eqref{cso} and \eqref{tedfu} are the same, which justifies (i).

To verify assertion (ii), note that for all $w\in K_\ph(\ox,\ov)$ the optimal value of the linear program
\begin{equation}\label{lp}
\max_{\lambda\in\Lambda(\ox,\ov)}\big\langle\lambda,\nabla^2f(\ox)(w,w)\big\rangle
\end{equation}
is finite, and thus this problem admits an optimal solution. We know from standard theory of linear programming that the set of optimal solutions to a linear program is a face of its feasible solution set. Denote by $G_1,\ldots,G_l$ all the finitely many faces of the polyhedral convex set $\Lambda(\ox,\ov)$. Select $\lm_i\in G_i$ for each $i=1,\ldots,l$ and define $E:=\{\lm_i|\;i=1,\ldots,l\}$. Choose further a positive number $r$ such that $E\subset\Lambda(\ox,\ov)\,\cap\,r\B$ and observe the equalities
\begin{equation}\label{gb01}
\max_{\lambda\in\Lambda(\ox,\ov)}\{\left\langle\lambda,\nabla^2f(\ox)(w,w)\right\rangle\}=\max_{\lambda\in E}\{\left\langle\lambda,\nabla^2f(\ox)(w,w)\right\rangle\}=\max_{\lambda\in[\Lambda(\ox,\ov)\,\cap\,r\B]}\{\left\langle\lambda,\nabla^2f(\ox)(w,w)\right\rangle\},
\end{equation}
which readily justify the second subderivative chain rule in (ii).
\end{proof}\vspace*{-0.05in}

{We finish this section by showing that the full subamenability of functions is preserved under various operations including the summation of functions and taking the second subderivative.}\vspace*{-0.05in}

\begin{Theorem}[\bf preservation of full subamenability under summation]\label{sumrule} Let $\ph:=\sum_{i=1}^{s}\ph_i$ on $\R^n$, where each $\ph_i\colon\R^n\to\oR$ is fully subamenable at $\ox\in\cap_{i=1}^s\dom\ph_i$, and let $\ov\in\ph(\ox)$. Impose the following subregularity qualification condition: there exist numbers $\kappa\in\R_+$ and $\ve>0$ such that
\begin{equation}\label{mq}
{\rm dist}\Big(x;\bigcap_{i=1}^{s}\dom\ph_i\Big)\le\kappa\,\sum_{i=1}^{s}{\rm dist}(x;\dom\ph_i)\quad\mbox{for all}\quad x\in\B_\ve(\ox).
\end{equation}
Then $\ph$ is fully subamenable at $\ox$ and its second subderivative is represented by
\begin{equation}\label{sumr3}
\d^2\ph(\ox,\ov)(w)=\max\Big\{\sum_{i=1}^s\d^2\ph_i(\ox,v_i)(w)\Big|\;v_i\in\sub\ph_i(\ox),\;\sum_{i=1}^s v_i=\ov\Big\},\quad w\in\R^n.
\end{equation}
\end{Theorem}\vspace*{-0.15in}
\begin{proof}
Since the functions $\ph_i$ are fully subamenable at $\ox$, for each $i\in\{1,\ldots,s\}$ there is a neighborhood $U_i$ of $\ox$ on which $\ph_i$ admits the representation $\ph_i=\vt_i\circ f_i$, where $\vt_i\colon\R^{p_i}\to\oR$ is convex piecewise linear-quadratic, where $f_i\colon\R^n\to\R^{p_i}$ is twice differentiable at $\ox$, and where the mapping $x\mapsto f_i(x)-\dom\vt_i$ is metrically subregular at $(\ox,0)$ with constant $\kappa_i\in\R_+$. Denote $p:=\sum_{i=1}^s p_i$ and define $\vt\colon\R^p\to\oR$ and $f\colon\R^n\to\R^p$ by, respectively,
$$
\vt(y_{p_1},\ldots,y_{p_s}):=\vt_1(y_{p_1})+\ldots+\vt_s(y_{p_s})\;\mbox{ as }\;y_{p_i}\in\R^{p_i}\;\mbox{ and }\;f(x):=\big(f_{1}(x),\ldots,f_s(x)\big).
$$
It is easy to check that the function $\vt$ is convex piecewise linear-quadratic and also that the mapping $x\mapsto f(x)-\dom\vt$ is metrically subregular at $(\ox,0)$ with constant $\bar{\kappa}:=\max\{\kappa\kappa_i|\;i=1,\ldots,s\}$ under the imposed subregularity qualification condition \eqref{mq}. This shows that $\ph=\vt\circ f$ on $U=\cap_{1=1}^s U_i$, and so $\ph$ is fully subamenable at $\ox$. The first-order chain rules from \eqref{dchain} and \eqref{1chain} applied to the composition $\ph=\vt\circ f$ allow us to arrive at the corresponding {\em sum rules}
\begin{equation}\label{sum-rule}
\sub\ph(\ox)=\sub\ph_1(\ox)+\ldots+\sub\ph_s(\ox)\quad\mbox{and}\quad\d\ph(\ox)(w)=\sum_{i=1}^s\d\ph_i(\ox)(w)\;\mbox{ for all }\;w\in\R^n.
\end{equation}

Next we verify representation \eqref{sumr3}. Pick $\ov\in\sub\ph(\ox)$ and choose subgradients $v_i\in\sub\ph_i(\ox)$ such that $\sum_{i=1}^s v_i=\ov$. It comes directly from the definition of the second subderivative that
$$
\d^2\ph(\ox,\ov)(w)\ge\sum_{i=1}^s\d^2\ph_i(\ox,v_i)(w)\;\mbox{ for all }\;w\in\R^n,
$$
which gives us the inequality ``$\ge$" in \eqref{sumr3}. To prove the opposite inequality, consider first the case where $w\notin K_\ph(\ox,\ov)$. It follows from definition \eqref{crit} and the subderivative sum rule in \eqref{sum-rule} that
$$
\sum_{i=1}^s\langle v_i,w\rangle=\langle\ov,w\rangle<\d\ph(\ox)(w)=\sum_{i=1}^s\d\ph_i(\ox)(w),
$$
which implies that there exists an index $i\in\{1,\ldots,s\}$ with $\langle v_i,w\rangle<\d\ph_i(\ox)(w)$. This ensures by the critical cone representation \eqref{domdf} that $\d^2\ph_i(\ox,v_i)(w)=\infty $. Since all $\ph_i$ are fully amenable at $\ox$, the second subderivatives $\d^2\ph_i(\ox,v_i)$ are proper and thus both sides in \eqref{sumr3} become $\infty$, which proves \eqref{sumr3} for all the noncritical vectors $w\notin K_\ph(\ox,\ov)$.

Consider now the remaining case where $w\in K_\ph(\ox,\ov)$. Let $\olambda\in\Lambda(\ox,\ov)$ be a vector realizing the maximum in the second subderivative representation \eqref{tedfu}. Take $\olambda_i\in\sub\vt_i(f_i(\ox))$ such that
$$
\olambda=(\olambda_1,\ldots,\olambda_s)\;\mbox{ and }\;\ov=\sum_{i=1}^s\ov_i\;\mbox{ with }\;\ov_i=\nabla f_i(\ox)^*\olambda_i.
$$\vspace*{-0.05in}
This implies by using \eqref{tedfu} the relationships
\begin{eqnarray*}
\d^2\ph(\ox,\ov)(w)&=&\d^2\vt\big(f(\ox),\olambda\big)\big(\nabla f(\ox)w\big)+\left\langle\olambda,\nabla^2f(\ox)(w,w)\right\rangle\\
&=&\sum_{i=1}^s\d^2\vt_i\big(f_i(\ox),\olambda_i\big)\big(\nabla f_i(\ox)w\big)+ \left\langle\olambda_i,\nabla^2f_i(\ox)(w,w)\right\rangle\\
&\le&\sum_{i=1}^s\d^2\ph_i(\ox,\ov_i)(w),
\end{eqnarray*}
which yield the inequality ``$\le$" in \eqref{sumr3} and thus complete the proof.
\end{proof}\vspace*{-0.05in}

{ The subregularity condition \eqref{mq} has been used in first-order theory of variational analysis under the name of ``metric qualification condition as, e.g., in \cite{io}. When in addition the sets in \eqref{mq} are convex and $x$ is taken from the entire space under consideration, \eqref{mq} reduces to ``linear regularity" that has been broadly studied in convex analysis and optimization. We refer the reader to \cite{bb,bd} for more details on the linear regularity in the convex framework.}

The final result of this section shows that the second subderivative of a fully subamenable function is fully subamenable itself on the whole space.\vspace*{-0.05in}

\begin{Corollary}[\bf full subamenability of second subderivatives]\label{fulsuba} Let $\ph\colon\R^n\to\oR$ be a fully subamenable function at $\ox$, and let $\ov\in\sub\ph(\ox)$. Then the second subderivative $w\mapsto\d^2\ph(\ox,\ov)(w)$ is a fully subamenable function at every $w\in\R^n$.
\end{Corollary}\vspace*{-0.15in}
\begin{proof} It follows from the proof of Corollary~\ref{cso0}(ii) that whenever $\olambda\in\Lambda(\ox,\ov)$ we have
\begin{equation}\label{fulsuba1}
\d^2\ph(\ox,\ov)(w)=\d^2\vt\big(f(\ox),\olambda)\big(\nabla f(\ox)w\big)+\max_{\lambda\in E}\big\langle\lambda,\nabla^2f(\ox)(w,w)\big\rangle,
\end{equation}
where $E:=\{\lambda_i|\:i=1,\ldots,p\}$ for some $\lambda_i\in\Lambda(\ox,\ov)$ and $p\in\N$. Define the function $\psi\colon\R^p\to\R$ by $\psi(y_1,\ldots,y_p):=\max\{y_1,\ldots,y_p\}$ and the mapping $g\colon\R^n\to\R^p$ by
$$
g(w):=\big(\big\langle\lambda_1,\nabla^2 f(\ox)(w,w)\big\rangle,\ldots,\big\langle\lambda_p,\nabla^2f(\ox)(w,w)\big\rangle\big).
$$
It is obvious that $\psi$ is convex piecewise linear-quadratic with $\dom\psi=\R^p$, that $g$ is a ${\cal C}^2$-smooth, and that MSQC \eqref{mscq} is satisfied at any point $w\in\R^n$ for the composition $\psi\circ g$. Thus the latter is fully subamenable at $w$. The second-order chain rule \eqref{fulsuba1} can be rewritten as
$$
\d^2\ph(\ox,\ov)(w)=\theta(w)+(\psi\circ g)(w)\;\mbox{ with }\;\theta(w):=\d^2\vt\big(f(\ox),\olambda\big)\big(\nabla f(\ox)w\big)\;\mbox{ for all }\;w\in\R^n.
$$
Remember that  both $\theta$ and $\psi\circ g$ are subamenable at any point $w\in \R^n$. To conclude that $\d^2\ph(\ox,\ov)$ is fully subamenable at $w$ by using Theorem~\ref{sumrule}, it remains to check that the subregularity qualification condition \eqref{mq} holds in our setting. This follows directly from the Hoffman lemma since the domains of $\theta$ and $\psi\circ g$ are polyhedral.
\end{proof}\vspace*{-0.25in}

\section{Second-Order Optimality Conditions for Composite Problems}\label{2opt}
\sce

Having in hand the developed calculus rules for second subderivatives, we are now in a position to derive {\em no-gap} second-order optimality conditions (i.e., such conditions where the difference between necessary and sufficient ones is in the replacement of the nonstrict inequality for its strict counterpart) for {\em composite optimization} problems written in the unconstrained format
\begin{equation}\label{op2}
\mbox{minimize }\;\ph_0(x)+(\vt\circ f)(x)\;\mbox{ over }\;x\in\R^n.
\end{equation}
Our assumptions here require that the mappings $\ph_0\colon\R^n\to\R$ and $f\colon\R^n\to\R^m$ are twice differentiable at $\ox$, and that $\vt\colon\R^m\to\oR$ is a convex piecewise linear-quadratic function. As already mentioned in Section~\ref{intro}, the possibility of taking the value $\vt(y)=\infty$ for the outer function in the composition from \eqref{op2} allows us to convert  constrained problems  into unconstrained problems. However, the realization of this approach to constrained optimization requires adequate generalized differential calculus to deal with extended-real-valued functions.

As discussed in Section~\ref{intro}, the composite format \eqref{op2} clearly covers classical problems of nonlinear programming, which correspond to the case where $\vt$ is the indicator function of a polyhedral convex  set. Another particular setting of \eqref{op2} is when $\vt\colon\R^m\to\oR$ is defined by \eqref{enlp}. As mentioned in Example~\ref{fmr}, this allows us to address  ENLPs the importance of which has been highly recognized in theoretical developments, computational methods, and applications dealing with broad areas of optimization including stochastic programming, robust optimization; see, e.g., Rockafellar and Wets \cite{rw} for further information. We also refer the reader to the more recent papers by Mordukhovich et al. \cite{mrs} and Do et al. \cite{dms}, devoted to the study of various stability issues, criticality of multipliers, and other aspects of ENLP important for numerical methods and applications.\vspace*{0.05in}

Now we derive {\em no-gap} second-order necessary and sufficient conditions for local optimality in composite problems described by fully subamenable functions. They are surely applied to ENLPs.\vspace*{-0.05in}

\begin{Theorem}[\bf no-gap second-order optimality conditions for fully subamenable composite problems]\label{nsop1}
Consider the composite optimization problem \eqref{op2}, where $\ph_0\colon\R^n\to\R$ and $f\colon\R^n\to\R^m$ are twice differentiable at $\ox$, and where $\vt\colon\R^m\to\oR$ is a convex piecewise linear-quadratic function with $f(\ox)\in\dom\vt$. Let $\ph:=\vt\circ f$, and let MSQC \eqref{mscq} hold at $\ox$ satisfying the stationary condition $0\in\nabla\ph_0(\ox)+\sub\ph(\ox)$. Take further $K_\ph(\ox,\ov)$ from \eqref{crit} with $\ov:=-\nabla\ph_0(\ox)$. The following hold:

{\bf(i)} If $\ox$ is a local minimizer of \eqref{op2}, then the second-order necessary condition
\begin{equation}\label{nopc1}
\d^2\vt\big(f(\ox),\bar\lm\big)\big(\nabla f(\ox)w\big)+\max_{\lambda\in\Lambda(\ox,\ov)}\big\langle\nabla^2_{xx} L(\ox,\lambda)w,w\big\rangle\ge 0\;
\mbox{ whenever }\;w\in K_\ph(\ox,\ov)
\end{equation}
is satisfied for any $\bar\lm\in\Lambda(\ox,\ov)$, where $L$ is the Lagrangian associated with \eqref{op2} and defined by $L(x,\lm):=\ph_0(x)+\la\lm,f(x)\ra-\vt^*(\lm)$ as $x\in\R^n$ and $\lm\in\R^m$.

{\bf(ii)} The validity of the second-order condition
\begin{equation}\label{sopc1}
\d^2\vt\big(f(\ox),\bar\lm\big)\big(\nabla f(\ox)w\big)+\max_{\lambda\in\Lambda(\ox,\ov)}\big\langle\nabla^2_{xx}L(\ox,\lambda)w,w \big\rangle>0 \mbox{ whenever } w\in K_\ph(\ox,\ov)\setminus\{0\}
\end{equation}
for any $\bar\lm\in\Lambda(\ox,\ov)$ amounts to the existence of numbers $\ell\ge 0$ and $\ve>0$ such that
\begin{equation}\label{quadg2}
\psi(x)\ge\psi(\ox)+\ell\|x-\ox\|^2\;\mbox{ if }\;x\in\B_{\ve}(\ox)
\end{equation}
with $\psi:=\ph_0+\vt\circ f$. In particular, condition \eqref{sopc1} is sufficient for local optimality of $\ox$ in \eqref{op2}.
\end{Theorem}\vspace*{-0.15in}
\begin{proof} Since $\ph_0$ is twice differentiable at $\ox$, it is easy to deduce from the definitions that
\begin{equation}\label{ns1}
\d^2(\ph_0+\ph)(\ox,0)(w)=\langle\nabla^2\ph_0(\ox)w,w\rangle+\d^2\ph(\ox,\ov)(w)\;\mbox{ for any }\;w\in\R^n.
\end{equation}
To verify (i), let $\ox$ be a local minimizer of \eqref{op2}, i.e., it provides a local minimum for $\psi=\ph_0+\ph$. It is an immediate consequence of definition \eqref{ssd} that $\d^2\psi(\ox,0)(w)\ge 0$ for all $w\in\R^n$. Applying the second-order sum rule \eqref{ns1} to $\psi$ and then the second-order chain rule \eqref{cso} to $\ph=\vt\circ f$, we arrive at the second-order necessary condition \eqref{nopc1} whenever $w\in K_\ph(\ox,\ov)$, which justifies (i).

To proceed with (ii), we first use Theorem~13.24(c) from Rockafellar and Wets \cite{rw} telling us that for any proper function $\psi\colon\R^n\to\oR$ the simultaneous fulfillment of the conditions $0\in\partial\psi(\ox)$ and $\d^2\psi(\ox,0)(w)>0$ when $w\ne 0$ amounts to the second-order growth condition  \eqref{quadg2}. Since $\psi=\ph_0+\vt\circ f$, we see from the elementary first-order subdifferential sum rule that the stationary point $\ox$ with $0\in\nabla\ph_0(\ox)+\partial(\vt\circ f)(\ox)$ satisfies $0\in\partial\psi(\ox)$. Furthermore, by Theorem~\ref{fsubam2} and the aforementioned sum rule for second subderivatives we get  $\dom\d^2\psi(\ox,0)=K_\ph(\ox,\ov)$. Applying as above the equality-type sum and chain rules for the second subderivative of $\psi$ shows that the condition $\d^2\psi(\ox,0)(w)>0$ amounts to \eqref{sopc1}. This completes the proof of the theorem.
\end{proof}\vspace*{-0.05in}

The no-gap second-order optimality conditions of Theorem~\ref{nsop1} are new while they can be derived in the same way by using the results in Rockafellar and Wets \cite{rw} (mainly the second-order chain rule in Theorem~13.14 therein) under  the metric regularity qualification condition \eqref{bqc}  as well as the ${\cal C}^2$-smoothness of $\ph_0$ and $f$ in \eqref{op2}. We also refer the reader to the concurrent preprint by Chieu et al. \cite{chnt}, where the second-order optimality conditions are obtained under a certain metric subregularity by a different approach for problems with the so-called ${\cal C}^2$-cone reducible constraints (in the sense of Bonnans and Shapiro \cite{bs}), which do not generally cover the case of subamenable compositions in \eqref{op2}. Note finally that the possibility of replacing ${\cal C}^2$-smoothness assumptions in conventional second-order optimality conditions for problems of nonlinear programming by merely twice differentiability of their data  has been already observed before; see \cite[Theorems~1.19, 1.20]{is}.\vspace*{0.05in}

{The following example demonstrates two simple composite problems for which the second-order optimality conditions
from Theorem~\ref{nsop1} can be applied while those in \cite[Theorem~13.14]{rw}, obtained under  the metric regularity qualification condition \eqref{bqc}, cannot be used.}\vspace*{-0.05in}
{
\begin{Example}[\bf second-order optimality conditions in the absence of metric regularity]\label{sof} {\rm Consider the composite optimization problem \eqref{op2} with the following initial data:

{\bf(a)} Take $f$ and $\vt$ from Example~\ref{fmr}(a) and choose $\ph_0(x)\equiv 0$ for $x\in \R^2$. It is easy to see that
$$
\vt(f(x))=\sm (x_1-x_2)^2 \quad\mbox{for all}\;\;x=(x_1,x_2)\in\R^2.
$$
This clearly indicates that $\ox=(0,0)$ is an optimal solution to \eqref{op2}. Remember that Example~\ref{fmr}(a) tells us that MSQC \eqref{mscq} is satisfied for the composite function $\ph=\vt\circ f$ at $\ox$, and thus this function is fully subamenable at $\ox$. Appealing to Theorem~\ref{nsop1}(i) tells us that the second-order necessary optimality condition \eqref{nopc1} holds for this problem. In fact, simple calculations show that
$$
K_{\ph}(\ox,\ov)=\big\{(w_1,w_2)\in\R^2\big|\;w_2=0\big\},\;\;\Lambda(\ox,\ov)=\big\{(\lm_1,\lm_2)\in\R^2\big|\;\lm_1= 0 \big\}\;\;\mbox{with}\;\;
\ov=(0,0).
$$
Moreover, for any $\lm=(\lm_1,\lm_2)\in\Lambda(\ox,\ov)$ and any $w=(w_1,w_2)\in\R^2$ we conclude that
$$
\big\langle\nabla^2_{xx}L(\ox,\lambda)w,w \big\rangle=0\quad\mbox{and}\quad\d^2\vt\big(f(\ox),\lm\big)(w)=w_1^2+\delta_{K_{\ph}(\ox,\ov)}(w),
$$
which readily confirms the fulfillment of \eqref{nopc1} for this composite problem. It is worth mentioning that any point $x=(x_1,x_2)\in\R^2$ with $x_1=x_2$ is also an optimal solution to \eqref{op2} at which the composite function $\ph$ is fully subamenable. Hence Theorem~\ref{nsop1}(i) can be utilized to achieve the necessary optimality condition for all such points while \cite[Theorem~13.14]{rw} cannot be used since the metric regularity qualification condition \eqref{bqc} fails at these points.

{\bf(b)} Take $f$ and $\vt$ from Example \ref{fmr}(b) and define the function $\ph_0\colon\R^3\to\R$ by $\ph_0(x):=-x_1+x^2_{2}+x^2_3$ with $x=(x_1,x_2,x_3)\in\R^3$. According to Example \ref{fmr}(b), MSCQ \eqref{mscq} holds at $\ox=(0,0,0)$ while the metric regularity qualification condition \eqref{bqc} fails at this point. Thus the function $\ph=\vt\circ f$ is fully subamenable (but not fully amenable) at $\ox$. It is easy to see that $\vt=\dd_{\R_-^3}$. Simple calculations show that $\ov=-\nabla\ph_0(\ox)=(1,0,0)$, that
\begin{eqnarray*}
K_{\ph}(\ox,\ov)&=&\big\{w=(w_1,w_2,w_3)\in\R^3\big|\;d\ph(\ox)(w)=w_1\big\}\\
&=&\big\{w=(w_1,w_2,w_3)\in\R^3\big|\;d\vt\big(f(\ox)\big)\big(\nabla f(\ox)w\big)=w_1\big\}\\
&=&\big\{w=(w_1,w_2,w_3)\in\R^3\big|\;\dd_{\R_-^3}(w_1,w_2,-w_3)=w_1\big\}=\{0\}\times\R_-\times\R_+,
\end{eqnarray*}
and that the Lagrange multiplier set is calculated by
$$
\Lambda(\ox,\ov)=\big\{(\lm_1,\lambda_2,\lm_3)\in\R_+^3\big|\;\lm_1+\lm_2-\lm_3=1\big\}.
$$
Since the latter set has the two extreme points $(1,0,0)$ and $(0,1,0)$, we get
\begin{eqnarray}
&&\max\big\{\big\la\nabla^2_{xx}L(\ox,\lambda)w,w\big\rangle\big|\;\lambda\in\Lambda(\ox,\ov)\big\}\nonumber\\
&=&\max\Big\{-\lm_3w_1^2+(2-\lm_1-\lm_3)w_2^2+(2-\lm_2-\lm_3)w_3^2\Big|\;\lm=(\lm_1,\lambda_2,\lm_3)\in\Lambda(\ox,\ov)\Big\}\label{lpma}\\
&=&\max\big\{w^2_2+2w^2_3,2w^2_2+w^2_3\big\}.\nonumber
\end{eqnarray}
Furthermore, it follows that $d^2\vt(f(\ox),\lm)=\delta_{\ss K_{\vt}(f(\ox),\lambda)}$ for any $\lm\in\Lambda(\ox,\ov)$. Thus the second-order sufficient condition \eqref{sopc1} holds at $\ox$, which ensures that the second-order growth condition \eqref{quadg2} is also satisfied at this point. It allows us to conclude that $\ox$ is a strict local minimizer of \eqref{op2} for this choice of the functions $\ph_0$, $\vt$, and $f$.}
\end{Example}\vspace*{-0.2in}

\color{black}\section{Calculating Proto-Derivatives of Subdifferentiable Mappings}\label{gra-der}
\sce\vspace*{-0.1in}

This section is devoted to deriving a precise calculus formula for the proto-derivative  of the subdifferential mapping associated with fully subamenable compositions. Recall that a set-valued mapping $S:\R^n\tto \R^m$ is {\em proto-differentiable} at $\ox$ for $\oy\in S(\ox)$ if $\gph S$ is geometrically derivable at $(\ox,\oy)$. When this condition holds for $S$, we refer to $DS(\ox,\oy)$ as the {\em proto-derivative} of $S$ at $\ox$ for $\oy$.

To achieve our goals, we implement a brilliant result first discovered by Rockafellar \cite{r90} for convex functions and then extended by Poliquin and Rockafellar \cite{pr96} to a significantly larger class of prox-regular and subdifferentially continuous functions. This result establishes a precise relationship between graphical derivatives of subdifferential mapping of a function and its second subderivative, which also leads us to new applications to parametric optimization developed in the next section.

To formulate this result, recall that a function $\psi\colon\R^n\to\oR$ is {\em prox-regular} at $\ox\in\dom\psi$ for $\ov\in\sub\psi(\ox)$ if it is l.s.c.\ around $\ox$ and there exist constants $\ve>0$ and $r>0$ such that for all $x\in\B_{\vep}(\ox)$ with $\psi(x)\le\psi(\ox)+\ve$ we have the condition
\begin{eqnarray*}
\psi(x)\ge\psi(u)+\la v,x-u\ra-\frac{r}{2}\|x-u\|^2\quad\mbox{whenever}\quad(u,v)\in(\gph\sub\psi)\cap\B_{\vep}(\ox,\ov).
\end{eqnarray*}
It is said that $\psi$ is {\em subdifferentially continuous} at $\ox$ for $\ov$ if the convergence $(x_k,v_k)\to(\ox,\ov)$ with $v_k\in\sub\psi(x_k)$ yields $\psi(x_k)\to\ph(\ox)$ as $k\to\infty$.

The aforementioned result by Poliquin and Rockafellar \cite[Theorem~6.1]{pr96} tells us that if $\psi\colon\R^n\to\oR$ is prox-regular and subdifferentially continuous at $\ox$ for $\ov\in\sub\psi(\ox)$, then $\psi$ is twice epi-differentiable at $\ox$ for $\ov$ if and only if $\sub\psi$ is proto-differentiable at $\ox$ for $\ov$ and its proto-derivative can be calculated by
\begin{equation}\label{gdfu}
\big(D\sub\psi\big)(\ox,\ov)(w)=\sm\sub\big(\d^2\psi(\ox,\ov)\big)(w)\;\mbox{ whenever }\;w\in\R^n.
\end{equation}

The fundamental relationship \eqref{gdfu} would allow us to employ the results obtained above for second subderivatives of fully subamenable functions if we show that such functions are prox-regular and subdifferentially continuous. This requires to assume in addition that the mapping $f$ in Definition~\ref{full} is ${\cal C}^2$-smooth around the reference point.\vspace*{-0.05in}

\begin{Proposition}[\bf prox-regularity and subdifferential continuity of fully subamenable functions]\label{sm} Let $\ph\colon\R^n\to\oR$ admit the representation $\ph=\vt\circ f$ locally around $\ox$, where $\vt\colon\R^m\to\oR$ is convex piecewise linear-quadratic, and where $f\colon\R^n\to\R^m$ is a ${\cal C}^2$-smooth mapping around $\ox$ with $f(\ox)\in\dom\vt$ under the fulfillment of MSQC \eqref{mscq} at $\ox$. Then $\ph$ is prox-regular and subdifferentially continuous at $\ox$ for any subgradient vector $\ov\in\sub\ph(\ox)$.
\end{Proposition}\vspace*{-0.15in}
\begin{proof} Since $f$ is ${\cal C}^2$-smooth around $\ox$, the subdifferential chain rule \eqref{1chain}  from Theorem~\ref{fchain} ensures the existence of $\ve>0$ such that
$$
\sub\ph(u)=\nabla f(u)^*\sub\vt\big(f(u)\big)\;\mbox{ whenever }\;u\in\B_\ve(\ox).
$$
Furthermore, the ${\cal C}^2$-smoothness of $f$ yields the boundedness property \eqref{subboundeq} for all $u\in\B_\ve(\ox)$ with the same constants on the right-hand side of the inequality in \eqref{subboundeq}. Picking now any $(u,v)\in(\gph\sub\ph)\cap\B_{\vep}(\ox,\ov)$ and appealing again to Corollary~\ref{subbound}, we find $\gg>0$ and $\lm\in \sub \ph(f(u))$ for which
\begin{equation*}
v=\nabla f(u)^*\lambda\quad\mbox{and}\quad\|\lambda\|\le\gg.
\end{equation*}
Combining the latter with the convexity of $\vt$ and the ${\cal C}^2$-smoothness of $f$ ensures the existence of $r>0$ such that
for any $(u,v)\in(\gph\sub\ph)\cap\B_{\vep}(\ox,\ov)$ and any $x\in\B_\ve(\ox)$ the relationships
\begin{eqnarray*}
\vt\big(f(x)\big)-\vt\big(f(u)\big)\ge\langle\lambda,f(x)-f(u)\rangle\ge\la\nabla f(u)^*\lambda,x-u\rangle-\frac{r}{2}\|x-u\|^2=\langle v,x-u\rangle-\frac{r}{2}\|x-u\|^2
\end{eqnarray*}
are satisfied. This verifies the prox-regularity of $\ph$ at $\ox$ for $\ov$. The subdifferential continuity of $\ph$ clearly follows from the fact that $\vt$ is continuous relative to its domain.
\end{proof}\vspace*{-0.05in}

Now we are ready to derive the following second-order chain rule for proto-derivatives of fully subamenable compositions.\vspace*{-0.05in}

\begin{Theorem}[\bf calculating proto-derivatives of fully subamenable functions]\label{gtosub} In the framework of Proposition~{\rm\ref{sm}}, assume that $(\ox,\ov)\in\gph\sub\ph$. Then the subgradient mapping $\sub\ph$ is proto-differentiable at $\ox$ for $\ov$ and whenever $w\in K_\ph(\ox,\ov)$ the proto-derivative of $\sub\ph$ is calculated by
\begin{equation}\label{gs}
\big(D\sub\ph\big)(\ox,\ov)(w)=\big\{\nabla^2\la\lambda,f\ra(\ox)w\big|\;\lambda\in\Lambda(\ox,\ov,w)\big\}+\nabla f(\ox)^*\big[\big(D\sub\vt\big)\big(f(\ox),\olambda\big)\big(\nabla f(\ox)w\big)\big],
\end{equation}
where $\Lambda(\ox,\ov,w)$ is the set of optimal solutions to the linear program \eqref{lp}, and where $\olambda$ is an arbitrary vector from $\Lambda(\ox,\ov)$. Furthermore, we have $\dom(D\sub\ph)(\ox,\ov)=K_\ph(\ox,\ov)$.
\end{Theorem}\vspace*{-0.15in}
\begin{proof} Theorem~\ref{sochain} tells us that the function $\ph$ is twice epi-differentiable at $\ox$ for $\ov$. As explained above, the latter amounts to saying that the mapping $\sub\ph$ is proto-differentiable at $\ox$ for $\ov$. Observe further from the second subderivative formula \eqref{cso2} that for all $w\in\R^n$ we have the representation $\d^2\ph(\ox,\ov)(w)=\psi_1(w)+\psi_2(w)$, where
\begin{equation*}
\psi_1(w):=\d^2\vt\big(f(\ox),\olambda\big)\big(\nabla f(\ox)w\big)\quad\mbox{and}\quad\psi_2(w):=\max_{\lambda\in[\Lambda(\ox,\ov)\,\cap\,r\B]}\big\langle\lambda,\nabla^2f(\ox)(w,w)\big\rangle
\end{equation*}
with $\bar\lm\in\Lambda(\ox,\ov)$ and $r>\bar r$, where $\bar r$ is taken from Corollary~\ref{cso0}(ii). Pick $w\in K_\ph(\ox,\ov)$, which ensures that $\nabla f(\ox)w\in K_\vt(f(\ox),\bar\lm)$. Since $\vt$ is convex piecewise linear-quadratic, so is $\d^2\vt(f(\ox),\olambda)$. Thus $\psi_1$ falls into the framework of Corollary~\ref{chpwlq} for which the subdifferential chain rule is always satisfied. This means that whenever $\nabla f(\ox)w\in K_\vt(f(\ox),\bar\lm)$ we get
\begin{equation}\label{gs2}
\sub\psi_1(w)=\nabla f(\ox)^*\sub\big[\d^2\vt\big(f(\ox),\olambda\big)\big]\big(\nabla f(\ox)w\big).
\end{equation}
Turing now to $\psi_2$, we conclude from Rockafellar and Wets \cite[Theorem~10.31]{rw} that this function is Lipschitz continuous and its subdifferential is calculated as
\begin{equation*}
\sub\psi_2(w)=\co\big\{2\nabla^2\la\lambda,f\ra(\ox)w\big|\;\lambda\in E_{\ss\rm{opt}}\big\},
\end{equation*}
where $E_{\ss\rm{opt}}$ is the set of optimal solutions to the problem
\begin{equation*}\label{gs3a}
\max_{\lambda\in[\Lambda(\ox,\ov)\,\cap\,r\B]}\big\langle\lambda,\nabla^2f(\ox)(w,w)\big\rangle.
\end{equation*}
Then \eqref{gb01} tells us that each optimal solution to this problem is an optimal solution to problem \eqref{du}. It results in  $E_{\ss\rm{opt}}=\Lambda(\ox ,\ov, w)\,\cap\,r\B$, where $\Lambda(\ox,\ov,w)$ is the set of optimal solutions to the linear program \eqref{lp}. Since this holds for any $r>\bar r$, taking the union  over all $r>\bar r$ implies that $E_{\ss\rm{opt}}=\Lambda(\ox,\ov, w)$. This leads us to the equalities
\begin{equation}\label{gs3}
\sub\psi_2(w)=\co\big\{2\nabla^2\la\lambda,f\ra(\ox)w\big|\;\lambda\in\Lambda(\ox,\ov,w)\big\}=\big\{2\nabla^2\la\lambda,f\ra(\ox)w\big|\;\lambda\in
\Lambda(\ox,\ov,w)\big\},
\end{equation}
where the last one comes from the convexity of the set $\{2\nabla^2\la\lambda,f\ra(\ox)w|\;\lambda\in\Lambda(\ox,\ov,w)\}$.
Employing the subdifferential sum rule and combining \eqref{gdfu}, \eqref{gs2}, and \eqref{gs3} imply that
\begin{eqnarray*}
\big(D\sub\ph\big)(\ox,\ov)(w)&=&\sm\sub\big[\d^2\ph(\ox,\ov)\big]=\sm\sub\psi_1(w)+\sm\sub\psi_2(w)\\
&=&\nabla f(\ox)^*\sub\big[\sm\d^2\vt\big(f(\ox),\olambda\big)\big](\nabla f(\ox)w)+\big\{\nabla^2\la\lambda,f\ra(\ox)w|\big|\;\lambda\in\Lambda(\ox,\ov,w)\big\}\\
&=&\nabla f(\ox)^*\big[\big(D\sub\vt\big)\big(f(\ox),\olambda\big)\big(\nabla f(\ox)w\big)\big]+\big\{\nabla^2\la\lambda,f\ra(\ox) w\big|\;\lambda\in\Lambda(\ox,\ov,w)\big\}
\end{eqnarray*}
for all $w\in\R^n$, where in the last equality we used \eqref{gdfu} for $\vt$.

It remains to verify the claimed formula for the domain of $(D\sub\ph)(\ox,\ov)$. Observe by Corollary~\ref{fulsuba} that
$d^2\ph(\ox,\ov)$ is a fully subamenable function. This allows us to deduce from the chain rule in \eqref{1chain} and Proposition~\ref{nes} that the subdifferential of the second subderivative function at any point of its domain is nonempty. Since $\dom d^2\ph(\ox,\ov)=K_\ph(\ox,\ov)$ due to \eqref{domdf}, we get $\dom(D\sub\ph)(\ox,\ov)=K_\ph(\ox,\ov)$ and thus complete the proof of the theorem.
\end{proof}\vspace*{-0.05in}

The proto-derivative formula for a fully subamenable function, obtained in \eqref{gs}, requires the proto-derivative of convex piecewise linear-quadratic functions. Now we provide a simple formula for the latter.\vspace*{-0.05in}

\begin{Proposition}[\bf proto-derivatives for convex piecewise linear-quadratic functions]\label{gdcp} Let $\vt\colon\R^m\to\oR$ be a convex piecewise linear-quadratic function with representation \eqref{PWLQ}, and let $(\oy,\ou)\in\gph\sub\vt$. Then
\begin{equation}\label{tan}
\big(D\sub\vt\big)(\oy,\ou)(w)=\bigcap_{i\in{\mathfrak J}(w)}\big\{y\in\R^m\big|\;y-A_i w\in N_{\ss{\cal K}_i(\oy,\ou_i)}(w)\big\}
\end{equation}
for all $w\in\R^m$, where ${\mathfrak J}(w):=\big\{i\in I(\oy)\big|\;w\in{\cal K}_i(\oy,\ou_i):= T_{\O_i}(\oy)\cap\{\ou_i\}^{\bot}\big\}$, where $\ou_i:=\ou-A_i\oy-a_i$, and where $I(\oy)$ is taken from \eqref{dfs}.
\end{Proposition}\vspace*{-0.15in}
\begin{proof} Remember that the second subderivative $d^2\vt(\oy,\ou)$ was calculated in \eqref{pwfor}, and that Corollary~\ref{fulsuba} tells us that it is a convex piecewise linear-quadratic function. Applying \eqref{gdfu} and the subdifferential formula for convex piecewise linear-quadratic functions from Rockafellar and Wets \cite[p.~487]{rw} verifies the claimed representation \eqref{tan}.
\end{proof}\vspace*{-0.05in}

If the function $\ph$ in Proposition~\ref{gdcp} is piecewise linear, i.e., $A_i=0$ in representation \eqref{PWLQ} for all $i=1,\ldots,s$, then the proto-derivative formula \eqref{tan} can be significantly simplified. Indeed, in this case we get from Proposition~13.9 in Rockafellar and Wets \cite{rw} that $d^2\vt(\oy,\ou)=\dd_{K_\vt(\oy,\ou)}$. This together with \eqref{gdfu} yields the equalities
$$
\big(D\sub\vt\big)(\oy,\ou)(w)=\sub\big(\sm\d^2\vt(\oy,\ou)\big)(w)=N_{K_\vt(\oy,\ou)}(w),\quad w\in\R^m,
 $$
where $K_\vt(\oy,\ou)=\dom d^2\vt(\oy,\ou)$ comes from \eqref{pwfor}. It immediately brings us to the following assertion.\vspace*{-0.05in}

\begin{Corollary}[\bf proto-derivative for a subclass of fully subamenable functions]\label{gtpi} In the framework of Theorem~{\rm\ref{gtosub}}, assume that $(\ox,\ov)\in\gph\sub\ph$ and that $\vt$ is piecewise linear. Then $\sub\ph$ is proto-differentiable at $\ox$ for $\ov$ and for any $w\in\R^n$ we have
\begin{equation}\label{gs7}
\big(D\sub\ph\big)(\ox,\ov)(w)=\big\{\nabla^2\la\lambda,f\ra(\ox)w\big|\;\lambda\in\Lambda(\ox,\ov,w)\big\}+\nabla f(\ox)^*N_{K_\vt(f(\ox),\olambda)}\big(\nabla f(\ox)w\big),
\end{equation}
where $\Lambda(\ox,\ov,w)$ and $\olambda$ are taken from Theorem~{\rm\ref{gtosub}}.
\end{Corollary}\vspace*{-0.05in}

{The next example provides an implementation of \eqref{gs7} for the function $\ph$ from Example~\ref{fmr}(b).\vspace*{-0.05in}

\begin{Example}[\bf proto-derivative calculations]\label{exaproto}{\rm Take $f$ and $\vt$ from Example~\ref{fmr}(b) and then choose $\ox:=(0,0,0)$ and $\ov:=(1,0,0)$. As shown in Example~\ref{sof}(b), the composite function $\ph=\vt\circ f$ is fully subamenable at $\ox$ with
$$
K_{\ph}(\ox,\ov)=\{0\}\times\R_-\times\R_+\quad\mbox{and}\quad\Lambda (\ox,\ov)=\big\{(\lm_1,\lambda_2,\lm_3)\in\R_+^3\big|\;\lm_1+\lm_2-\lm_3=1\big\}.
$$
Since the Lagrange multiplier set $\Lambda(\ox,\ov)$ has the two extreme points $(1,0,0)$ and $(0,1,0)$, for any $w=(w_1,w_2,w_3)\in K_{\ph} (\ox,\ov)$ this set is calculated by
$$
\Lambda(\ox,\ov,w)=
\begin{cases}
\{(1,0,0)\}&\mbox{if}\;|w_2|>w_3,\\
\{(0,1,0)\}&\mbox{if}\;|w_2|<w_3,\\
\big\{(\lm_1,\lm_2,0)\big|\;\lm_1+\lm_2=1,\,\lm_1\ge0,\,\lm_2\ge 0\}&\mbox{if}\;|w_2|=w_3,\\
\end{cases}
$$
where $\Lambda(\ox,\ov,w)$ is the set of optimal solution to \eqref{lpma}. Choose $\olm:=(1,0,0)$ and remember from Example~\ref{sof}(b) that $\vt$ can be equivalently written as $\vt=\dd_{\R_-^3}$, and so that $\vt$ is piecewise linear. Thus
$$
K_{\vt}(f(\ox),\olambda)=\big\{u=(u_1,u_2,u_3)\in\R^3\big|\;\dd_{\R_-^3}(u)=u_1\big\}=\{0\}\times\R_{-}\times\R_{-}.
$$
Employing now \eqref{gs7}, we conclude for any $w=(w_1,w_2,w_3)\in \R^3$ that
\begin{eqnarray*}
&&\big(D\sub\ph\big)(\ox,\ov)(w)=\\
&&
\begin{cases}
\big\{(0,0,-\beta)\big|\;\beta\le 0\big\}&\mbox{if}\;w_1=0,\,|w_2|>w_3=0,\\
\{(0,0,0)\}&\mbox{if}\;w_1=0,\,|w_2|>w_3>0,\\
(0,0,-w_3)+\big\{(0,\al,0)\big|\;\al\ge 0\big\}&\mbox{if}\;w_1=0,\,w_3>|w_2|=0,\\
\{(0,0,-w_3)\}&\mbox{if}\;w_1=0,\,w_3>|w_2|>0,\\
\big\{(0,0,-\beta w_3)\big|\;\beta\in[0,1]\big\}+\big\{(0,\al,-\beta)\big|\;\al\ge 0,\,\beta\le 0\big\}&\mbox{if}\;w_1 =0,\,|w_2|=w_3=0,\\
\big\{(0,0,-\beta w_3)\big|\;\beta\in[0,1]\big\}&\mbox{if}\;w_1=0,\,|w_2|=w_3>0,\\
\emp&\mbox{if}\;w\notin K_{\ph}(\ox,\ov).
\end{cases}
\end{eqnarray*}}
\end{Example}}\vspace*{-0.05in}

We end this section by providing some comments about the obtained proto-derivative formulas \eqref{gs} and \eqref{gs7}. For fully amenable functions (i.e., under the metric regularity qualification condition) the chain rule in \eqref{gs} was first obtained by Poliquin and Rockafellar \cite[Proposition~2.10]{pr93}. Their result did not draw much attention at that time, but the recent progress in parametric optimization has revived the importance of finding proto-derivatives for important classes of constrained optimization problems. The new effort to calculate the graphical derivatives of subdifferential mappings, which is a weaker property than proto-derivatives, under metric subregularity was undertaken by Gfrerer and Outrata \cite{go16} in the framework of nonlinear programming. The approach therein is very different from ours and did not achieve the proto-differentiability of subdifferential mappings as our result in Theorem~\ref{gs}.
In fact, we obtained the proto-derivatives formulas for fully subamenable functions as a by product of their second subderivative formulas and their twice epi-differentiability that were established under MSQC \eqref{mscq}.\vspace*{-0.15in}

\section{Applications to Parametric Optimization}\label{param}
\sce\vspace*{-0.1in}

In this section we provide applications of the second-order developments presented above to some important topics in parametric optimization. These topics mainly concern the uniqueness of Lagrange multipliers and stability properties for KKT systems in composite optimization, which play a prominent role in the design and justification of numerical algorithms.

We begin with a complete characterization of {\em uniqueness of Lagrange multipliers} for a general class of composite optimization problems \eqref{op2}, which is an extension of the recent result by Mordukhovich and Sarabi \cite[Theorem~3.1]{ms19} obtained for constrained optimization problems. The first-order optimality conditions
for the composite problem \eqref{op2} are given by
\begin{equation}\label{foc}
\nabla_x L(x,\lm)=0,\;\;\lm\in\sub\vt\big(f(x)\big),\quad x\in\R^n,
\end{equation}
via the Lagrangian $L(x,\lm)=\ph_0(x)+\la\lm,f(x)\ra-\vt^*(\lm)$
{and  can be equivalently expressed as $(0,0)\in (\sub_x L(\ox,\olm),\sub_\lm(-L)(\ox,\olm))$.}
The optimality conditions \eqref{foc} motivate us to consider the mapping $G\colon\R^n\times\R^m\tto\R^n\times\R^m$ defined by
\begin{equation}\label{GKKT}
G(x,\lm):=\left[\begin{array}{c}
\nabla_x L(x,\lm)\\
-f(x)
\end{array}
\right]+\left[\begin{array}{c}
0\\
(\sub\vt)^{-1}(\lm)
\end{array}
\right].
\end{equation}

Recall that a set-valued mapping $S\colon\R^n\tto\R^m$ is {\em strongly metrically subregular} at $(\ox,\oy)\in\gph S$ if there exist a constant $\kappa\in\R_+$ and a neighborhood $U$ of $\ox$ such that the estimate
\begin{equation}\label{strong-mr}
\|x-\ox\|\le\kappa\,{\rm dist}\big(\oy;S(x)\big)\quad\mbox{for all}\;\;x\in U
\end{equation}
is satisfied. The {\em Levy-Rockafellar criterion} (see, e.g., \cite[Theorem~4E.1]{dr} and the commentaries therein) tells us that the set-valued mapping $S$ is strongly metrically subregular at $(\ox,\oy)$ if and only if the following implication in terms of the graphical derivative \eqref{gder} holds:
\begin{equation}\label{sms8}
0\in DS(\ox,\oy)(w)\implies w=0.
\end{equation}

We use this criterion in the next theorem to characterize the uniqueness of Lagrange multipliers for a general class of problems in composite optimization.\vspace*{-0.05in}

\begin{Theorem}[\bf characterization of the uniqueness of Lagrange multipliers]\label{unique} Consider the composite optimization problem \eqref{op2}, where $\ph_0\colon\R^n\to\R$ and $f\colon\R^n\to\R^m$ are twice differentiable at $\ox$, and where $\vt\colon\R^m\to\oR$ is convex piecewise linear-quadratic with $f(\ox)\in\dom\vt$. Let $(\ox,\bar\lm)$ be a solution to the KKT system \eqref{foc}. Then the following conditions are equivalent:

{\bf(i)} The set of Lagrange multipliers \eqref{lagn} with $\ov:=-\nabla\ph_0(\ox)$ is a singleton $\Lambda(\ox,\ov)=\{\bar\lm\}$.

{\bf(ii)} The qualification condition
\begin{equation}\label{gf02}
\big(D\sub\vt\big)\big(f(\ox),\bar\lm\big)(0)\cap\ker\nabla f(\ox)^*=\{0\}
\end{equation}
holds with $(D\sub\vt)(f(\ox),\bar\lm)(0)=K_\vt(f(\ox),\bar\lm)^*$, where $K_\vt(f(\ox),\bar\lm)$ is taken from \eqref{crit}.
\end{Theorem}\vspace*{-0.15in}
\begin{proof} Using $G$ from \eqref{GKKT}, define the set-valued mapping $G_{\ox}\colon\R^m\tto\R^n\times\R^m$ by $G_{\ox}(\lm):=G(\ox,\lm)$. Since $(\ox,\bar\lm)$ is a solution to \eqref{foc}, we have $(\bar\lm,(0,0))\in\gph G_{\ox}$. Let us first establish two claims about subregularity behavior of $G_{\ox}$. The first claim shows that this mapping is always metrically subregular at the reference point.\\[0.07ex]
{\bf Claim A:} {\em The set-valued mapping $G_{\ox}$ is metrically subregular at $(\bar\lm,(0,0))$.}

To verify it, recall that $\olm\in\sub\vt(f(\ox))$ and deduce from the proof Theorem~11.14(b) in Rockafellar and Wets \cite{rw} that the set $\gph\partial\vt$ is a union of finitely many polyhedral convex sets. Then it follows from the seminal result by Robinson \cite{rob} that the mapping $(\sub\vt)^{-1}$ is metrically subregular at $(\olm,f(\ox))$, which ensures the existence of a constant $\kappa\in\R_+$ and a neighborhood $U$ of $\olm$ for which the distance estimate
$$
{\rm dist}\big(\lm;\sub\vt(f(\ox))\big)\le\kappa\,{\rm dist}\big(f(\ox);(\sub \vt)^{-1}(\lm)\big)\quad\mbox{for all}\;\;\lm\in U
$$
holds. Observe that $G_{\ox}^{-1}(0,0)=\Lambda(\ox,\ov)$ and that $\sub\vt(f(\ox))$ is a polyhedral convex set. This tells us that $\Lambda(\ox,\ov)$ is the intersection of two polyhedral convex sets. Employing again the Hoffman lemma confirms that there exists a constant $\ell\in\R_+$ such that
\begin{eqnarray*}
{\rm dist}\big(\lm;G_{\ox}^{-1}(0,0)\big)={\rm dist}\big(\lm;\Lambda(\ox,\ov)\big)&\le&\ell\,\Big(\|\nabla_x L(\ox,\lm)\|+{\rm dist}\big(\lm;\sub\vt(f(\ox))\big)\Big)\\
&\le&\ell\,\Big(\|\nabla_x L(\ox,\lm)\|+\kappa\,{\rm dist}\big(f(\ox);(\sub\vt)^{-1}(\lm)\big)\Big)\\
&\le&\max\{\ell,\ell\kappa\}\,{\rm dist}\big((0,0);G_{\ox}(\lm)\big)\quad\mbox{for all}\;\;\lm\in U.
\end{eqnarray*}
This yields the metric subregularity of $G_{\ox}$ and thus justifies Claim~A.\vspace*{0.05in}

The next claim characterizes the strong metric subregularity of $G_{\ox}$ at the reference point.\\[0.05in]
{\bf Claim~B:} {\em The mapping $G_{\ox}$ is strongly metrically subregular at $(\bar\lm,(0,0))$ if and only if the qualification condition \eqref{gf02} is satisfied.}

To verify this claim, observe by the direct calculation that
$$
DG_{\ox}\big(\olm,(0,0)\big)(u)=\left[\begin{array}{c}
\nabla f(\ox)^*u\\
0
\end{array}
\right]+\left[\begin{array}{c}
0\\
D(\sub\vt)^{-1}\big(\olm,f(\ox)\big)(u)
\end{array}
\right]\quad\mbox{for all}\;\;u\in\R^m.
$$
Combining this with the Levy-Rockafellar criterion \eqref{sms8} justifies Claim~B.\vspace*{0.05in}

Now we continue with the proof of the theorem. To verify implication (i)$\implies$(ii), note that $G_{\ox}^{-1}(0,0)=\Lambda(\ox,\ov)=\{\bar\lm\}$. It follows from Claim~A that the mapping $G_{\ox}$ is strongly metrically subregular at $(\bar\lm,(0,0))$. Appealing now to Claim~B shows that condition \eqref{gf02} is satisfied. Furthermore, it follows from \eqref{tan} and representation \eqref{PWLQ} that
\begin{eqnarray*}
\big(D\sub\vt\big)\big(f(\ox),\bar\lm\big)(0)&=&\bigcap_{i\in{\mathfrak J}(0)}\big\{y\in\R^m\big|\;y\in N_{\ss{\cal K}_i(f(\ox),\olm_i)}(0)\big\}\\
&=&\bigcap_{i\in{\mathfrak J}(0)}N_{\ss{\cal K}_i(f(\ox),\olm_i)}(0)\\
&=&N_{\ss\cup_{i\in{\mathfrak J}(0)}{\cal K}_i(f(\ox),\olm_i)}(0)=N_{K_\vt(f(\ox),\bar\lm)}(0)=K_\vt\big(f(\ox),\bar\lm\big)^*,
\end{eqnarray*}
where ${\mathfrak J}(0):=\{i\in I(f(\ox))|\;0\in{\cal K}_i(f(\ox),\lm_i):=T_{\O_i}(f(\ox))\cap\{\olm_i\}^{\bot}\}$, where $\olm_i:=\olm-A_i f(\ox)-a_i$, and where $I(f(\ox))$ is taken from \eqref{dfs}. This verifies (ii).

It remains to check that implication (ii)$\implies$(i) holds. Indeed, it follows from Claim~B that the mapping $G_{\ox}$ is strongly metrically subregular at $(\bar\lm,(0,0))$. This yields the existence of a neighborhood $U$ of $\olm$ for which we have the equalities
$$
\Lambda(\ox,\ov)\cap U=G_{\ox}^{-1}(0,0)\cap U=\{\bar\lm\}.
$$
Using this and the convexity of the Lagrange multiplier set $\Lambda(\ox,\ov)$ tells us that $\Lambda(\ox,\ov)=\{\bar\lm\}$, which verifies (i) and thus completes the proof of the theorem.
\end{proof}\vspace*{-0.05in}

Looking at the qualification condition \eqref{gf02}, observe by \eqref{dfs} that
$$
K_\vt\big(f(\ox),\bar\lm\big)=\big\{u\in\R^m\big|\;\d\vt\big(f(\ox)\big)(u)=\langle u,\bar\lm\rangle\big\}\subset\dom\d\vt\big(f(\ox)\big)=T_{\ss\dom\vt}\big(f(\ox)\big),
$$
which brings us to the inclusion $N_{\ss\dom\vt}(f(\ox))\subset K_\vt(f(\ox),\bar\lm)^*$. This tells us that \eqref{gf02} yields the constraint qualification \eqref{gf07} of the Robinson type. Consequently, we  refer to condition \eqref{gf02} in what follows as to the {\em strong Robinson constraint qualification}.

Next we proceed with a characterization of strong metric subregularity of the mapping $G$ from \eqref{GKKT}. The equivalence between (i) and (ii) in the following theorem has been recently proved in Burke and Engle \cite[Theorem~5.1]{be} by a different method. {Such a result was first established by Dontchev and Rockafellar \cite[Theorem~2.6]{dr97} for classical nonlinear programs and then was extended by Ding et al.\cite{dsz} to ${\cal C}^2$-reducible constrained optimization problems; see also \cite{dms} for similar results obtained for ENLPs.} Our approach is based on the fundamental relationship \eqref{gdfu} and the characterization of the uniqueness of Lagrange multipliers established in Theorem~\ref{unique}. It allows us to provide a simpler proof of (i)$\iff$(ii) and add the new equivalence (iii).\vspace*{-0.05in}

\begin{Theorem}[\bf characterizations of strong metric subregularity of KKT systems]\label{isos}
Given the composite optimization problem \eqref{op2}, assume that $\ph_0\colon\R^n\to\R$ and $f\colon\R^n\to\R^m$ are ${\cal C}^2$-smooth around $\ox$, and that $\vt$ is convex piecewise linear-quadratic with $f(\ox)\in\dom\vt$. Let $(\ox,\olm)$ be a solution to the KKT system \eqref{foc}. Then the following conditions are equivalent:

{\bf(i)} The mapping $G$ from \eqref{GKKT} is strongly metrically subregular at $((\ox,\olm),(0,0))$, and $\ox$ is a local minimizer for \eqref{op2}.

{\bf(ii)} $\Lambda(\ox,\ov)=\{\olm\}$ with $\ov=-\nabla f_0(\ox)$, and the second-order sufficient condition
\begin{equation}\label{sscc}
 \d^2\vt\big(f(\ox),\bar\lm\big)\big(\nabla f(\ox)w\big)+\langle\nabla^2_{xx}L(\ox,\olm)w,w\rangle>0
\end{equation}
holds for all $w\in\R^n\setminus\{0\}$ satisfying $\nabla f(\ox)w\in K_\vt\big(f(\ox),\olm\big)$.

{\bf(iii)} The second-order sufficient condition \eqref{sscc} and the strong Robinson constraint qualification \eqref{gf02} are both satisfied.
\end{Theorem}\vspace*{-0.15in}
\begin{proof} Let us first verify the following claim of its own interest.\\[0.03in]
{\bf Claim:} {\em The set-valued  mapping $G$ defined in \eqref{GKKT} is strongly metrically subregular at $((\ox,\olm),(0,0))$ if and only if we have the implication
\begin{equation}\label{isos1}
\begin{cases}
\nabla^2_{xx}L(\ox,\olm)w+\nabla f(\ox)^*u=0,\\
u\in\big(D\sub\vt\big)\big(f(\ox),\olm\big)\big(\nabla f(\ox)w\big)
\end{cases}
\implies w=0,\;u=0.
\end{equation}}

To check it, observe that the graphical derivative of $G$ is calculated by
$$
DG\big((\ox,\olm),(0,0)\big)(w,u)=\left[\begin{array}{c}
\nabla^2_{xx}L(\ox,\olm)w+\nabla f(\ox)^*u\\
-\nabla f(\ox)w
\end{array}\right]
+\left[\begin{array}{c}
0\\\big(D\sub\vt\big)^{-1}\big(\olm,f(\ox)\big)(u)\\
\end{array}
\right].
$$
Combining the latter with \eqref{sms8} readily justifies this claim.

The equivalence between (ii) and (iii) is a direct consequence of Theorem~\ref{unique}. We now proceed with verifying implication (i)$\implies$(iii). Observe that (i) yields \eqref{gf02}. Indeed, if $u\in (D\sub\vt)(f(\ox),\bar\lm)(0)\cap\ker\nabla f(\ox)^*$, then we have by \eqref{isos1} that $u=0$ and thus get \eqref{gf02}. This tells us by Proposition~\ref{unique} that $\Lambda(\ox,\ov)=\{\olm\}$. As discussed above, \eqref{gf02} implies the constraint qualification \eqref{gf07}, which ensures by Proposition~\ref{equi} that MSQC \eqref{mscq} holds. Since $\ox$ is a local minimizer of \eqref{op2}, we deduce from Theorem~\ref{nsop1}(i), $\Lambda(\ox,\ov)=\{\olm\}$, and the equivalence in \eqref{pj98} that the second-order necessary condition
$$
d^2\vt\big(f(\ox),\bar\lm\big)\big(\nabla f(\ox)w\big)+\langle\nabla^2_{xx}L(\ox,\olm)w,w\rangle\ge 0
$$
is fulfilled for all $w\in\R^n$ with $\nabla f(\ox)w\in K_\vt(f(\ox),\olm)$.

To verify (iii), it remains to show that for $w\ne 0$ the above inequality is strict. Arguing by contradiction. suppose that there exists $\ow\ne 0$ with $\nabla f(\ox)\ow\in K_\vt(f(\ox),\olm)$ and
$$
d^2\vt\big(f(\ox),\bar\lm\big)\big(\nabla f(\ox)\ow\big)+\langle\nabla^2_{xx}L(\ox,\olm)\ow,\ow\rangle=0.
 $$
This implies that the vector $\ow$ is a local minimizer for the problem
\begin{equation*}
\min_{w\in\R^n}\sm\big[\la\nabla_{xx}^2L(\ox,\olm)w,w\ra+d^2\vt\big(f(\ox),\bar\lm\big)\big(\nabla f(\ox)w\big)\big].
\end{equation*}
Applying the subdifferential Fermat rule to the cost function above and then using the subdifferential chain rule from Corollary~\ref{chpwlq} tell us that
\begin{eqnarray*}
0&\in&\nabla_{xx}^2L(\ox,\olm)\ow+\nabla f(\ox)^*\sub\big[\sm\d^2\vt\big(f(\ox),\bar\lm\big)\big]\big(\nabla f(\ox)\ow\big)\\
&=&\nabla_{xx}^2L(\ox,\olm)\ow+\nabla f(\ox)^*\big(D\sub\vt\big)\big(f(\ox),\bar\lm\big)\big(\nabla f(\ox)\ow\big),
\end{eqnarray*}
where the last equality comes from \eqref{gdfu}. This together with \eqref{isos1} gives us $\ow= 0$, which is a contradiction that verifies (iii).

To justify the final implication (iii)$\implies$(i), we get from Theorem~\ref{unique} that $\Lambda(\ox,\ov)=\{\olm\}$. Taking it into account along with the second-order sufficient optimality condition \eqref{sscc} and then employing Theorem~\ref{nsop1}(ii) tell us that $\ox$ is a local minimizer of \eqref{op2}. We show now that $G$ from \eqref{GKKT} is strongly metrically subregular at $((\ox,\olm),(0,0))$. To proceed, pick a pair $(w,u)\in\R^n\times\R^m$ satisfying on the left-hand side of implication \eqref{isos1}. This brings us to
\begin{equation}\label{ip01}
\langle\nabla^2_{xx}L(\ox,\olm)w,w\rangle+\la u,\nabla f(\ox)w\ra=0\quad\mbox{and}\quad u\in\big(D\sub\vt\big)\big(f(\ox),\olm\big)\big(\nabla f(\ox)w\big).
\end{equation}
The second relationship together with \eqref{gdfu} yields $u\in\sub\big[\sm d^2\vt(f(\ox),\bar\lm)\big](\nabla f(\ox)w)$, which in turn implies that $\nabla f(\ox)w\in\dom\d^2\vt(f(\ox),\bar\lm)=K_\vt(f(\ox),\olm)$. Moreover, the convexity of the function $v\mapsto\d^2\vt(f(\ox),\bar\lm)(v)$ ensures by the subdifferential construction of convex analysis that
$$
\big\la u,v-\nabla f(\ox)w\ra\le\sm\d^2\vt\big(f(\ox),\bar\lm\big)(v)-\sm d^2\vt\big(f(\ox),\bar\lm\big)\big(\nabla f(\ox)w\big)\;\mbox{ for all }\;v\in\R^m.
$$
Pick further $\ve\in(0,1)$ and take $v:=(1\pm\ve)\nabla f(\ox)w$. Since the second subderivative is positive homogeneous of degree $2$,
we arrive at
$$
\pm\la u,\nabla f(\ox)w\ra\le\frac{\ve\pm 2}{2}\d^2\vt\big(f(\ox),\bar\lm\big)\big(\nabla f(\ox)w\big).
$$
Letting $\ve\dn 0$ clearly gives us the equality $d^2\vt(f(\ox),\bar\lm)(\nabla f(\ox)w)=\la u,\nabla f(\ox)w\ra$. Combining this and \eqref{ip01} tells us that
$$
\big\langle\nabla^2_{xx}L(\ox,\olm)w,w\big\rangle+d^2\vt\big(f(\ox),\bar\lm\big)\big(\nabla f(\ox)w\big)=0\;\mbox{ with }\;\nabla f(\ox)w\in K_\vt\big(f(\ox),\olm\big),
$$
which results in $w=0$ due to the second-order condition \eqref{sscc}. Letting $w=0$ in \eqref{ip01} yields
$$
u\in\big(D\sub\vt\big)\big(f(\ox),\bar\lm\big)(0)\cap\ker\nabla f(\ox)^*.
$$
Combined with \eqref{gf02}, the latter inclusion implies that $u=0$, and thus implication \eqref{isos1} holds. Thus the mapping $G$ from \eqref{GKKT} is strongly metrically subregular at $((\ox,\olm),(0,0))$, which verifies (i) and therefore completes the proof of the theorem.
\end{proof}\vspace*{-0.2in}

\section{Concluding Remarks}\label{conc}
\sce\vspace*{-0.1in}

This paper provides a comprehensive variational study of a major class of composite optimization problems. The underlying theme of our study is a systematic usage of a novel version of metric subregularity conditions in contrast to the conventional metric regularity ones. The developed approach allows us to obtain new results of first-order and second-order variational analysis with applications to optimization, and also to improve with simplified proofs important achievements in the case of metric regularity. A major class of compositions introduced and investigated in this paper consists of fully subamenable functions while being based on metric subregularity. For this class we develop second-order variational analysis with applications to optimization and stability at the same level of perfection as previously known for fully amenable one, which strongly depends on metric regularity.

Our future research, which is partly implemented in \cite{mms}, aims at overcoming a polyhedral structure of fully subamenable compositions with replacing it by the concept of parabolic regularity, which encompasses fairly general nonpolyhedral settings. In this way we indent to develop comprehensive second-order calculus rules with broad applications to conic programming, various stability issues, as well as to the design and justification of numerical algorithms in optimization.\vspace*{0.05in}

{\bf Acknowledgements.} Research of the first and second authors was partly supported by the USA National Science Foundation under grants DMS-1512846 and DMS-1808978, and by the USA Air Force Office of Scientific Research under grant \#15RT0462. The research of the second author was also supported by the Australian Research Council under Discovery Project DP-190100555.

\end{document}